\documentclass{siamart250211}
\usepackage{graphicx} 
\usepackage{amsmath,amstext,amssymb,bm}
\usepackage{leftidx}

\usepackage{soul} 
\usepackage{tikz}
\usetikzlibrary{shapes,arrows}
\usepackage{mathrsfs}
\usepackage{epstopdf}

\usepackage{multirow}
\usepackage{tabularx}
\usepackage[shortlabels]{enumitem}
\usepackage[config, labelfont={bf}]{caption,subfig}
\DeclareGraphicsExtensions{.jpg,.pdf,.png,.jpeg, .eps}

\usepackage{xcolor}

\usepackage{hyperref}
\usepackage{cleveref}
\hypersetup{
	colorlinks   = true,     
	urlcolor     = blue,     
	linkcolor    = red,      
	citecolor    = red,      
	breaklinks   = true,     
	pdfborder    = {0 0 0}   
}

\numberwithin{equation}{section}
\newtheorem{remark}{Remark}[section]

\allowdisplaybreaks[4]

\newcommand{\cQ}{\mathcal{Q}}

\renewcommand{\div}{\mbox{\rm div\,}}

\newcommand{\cF}{\mathcal{F}}

\newcommand{\mH}{\mathbb{H}}
\newcommand{\mP}{\mathbb{P}}
\newcommand{\mE}{\mathbb{E}}
\newcommand{\mV}{\mathbb{V}}

\newcommand{\Ome}{\Omega}

\newcommand{\nab}{\nabla}
\newcommand{\vu}{{\bf u}}

\newcommand{\vH}{{\bf H}}

\newcommand{\vL}{{\bf L}}

\newcommand{\vv}{{\bf v}}
\newcommand{\ve}{{\bf e}_{\vu}}

\newcommand{\vw}{{\bf w}}
\newcommand{\vA}{{\bf A}}
\newcommand{\vG}{{\bf G}}

\newcommand{\pphi}{\pmb{\phi}}

\begin{document}
	
	\title{A mixed finite element method for the stochastic Boussinesq equations with multiplicative noise}
	\markboth{LIET VO}{Fully discrete FEM for stochastic Boussinesq system}

	\author{
		Liet Vo\thanks{
			School of Mathematical and Statistical Sciences, The University of Texas Rio Grande Valley, Edinburg, TX 78539, U.S.A. (liet.vo@utrgv.edu). This author was partially supported by the NSF grant DMS-2530211.}
	}
	
	\maketitle
	
	\begin{abstract}
		This work investigates a fully discrete mixed finite element method for the stochastic Boussinesq system driven by multiplicative noise. 
		The spatial discretization is performed using a standard mixed finite element method, while the temporal discretization is based on a semi-implicit Euler-Maruyama scheme. 
		By combining a localization technique with high-moment stability estimates, we establish error bounds for the velocity, pressure, and temperature approximations. 
		As a direct consequence, we prove convergence in probability for the fully discrete method in both $L^2$ and $H^1$-type norms. 
		Several numerical experiments are presented to validate the theoretical error estimates and demonstrate the effectiveness of the proposed scheme.
	\end{abstract}

	
	
	\section{Introduction}
	\label{sec-1}
	We consider  the  following stochastic Boussinesq system with multiplicative noise
	\begin{subequations}\label{eq1.1}
		\begin{alignat}{2} \label{eq1.1a}
			d\vu &=\bigl[\nu\Delta \vu - (\vu\cdot \nab)\vu -\nabla p + \theta {\bf e}_2 \bigr] dt +  \vG_1(\vu)d W_1(t)  &&\qquad\mbox{a.s. in}\, D_T,\\
			d\theta &=\bigl[\mu\Delta \theta - \vu \cdot\nab\theta \bigr] dt +  G_2(\theta)d W_2(t)  &&\qquad\mbox{a.s. in}\, D_T,\\
			\div \vu &=0 &&\qquad\mbox{a.s. in}\, D_T,\label{eq1.1b}\\
			\vu(0)&= \vu_0 &&\qquad\mbox{a.s. in}\, D,\label{eq1.1d}
		\end{alignat}
	\end{subequations}
	where $D = (0,L)^2 \subset \mathbb{R}^2 \,$ represents a period of the periodic domain in $\mathbb{R}^2$, $\vu$, $p$, and $\theta$ stand for respectively the velocity field, the pressure, and the temperature of the fluid. $\{W_{i}(t); t\geq 0\}_{i=1,2}$ denotes real-valued Wiener processes. The constant $\nu>0$ denotes the viscosity of the fluid, and $\mu>0$ represents the thermal diffusivity. In addition, $\vG_1$ and $G_2$ are the diffusion coefficients (see section \ref{sub5.2.2} for their precise definition).  Here we seek periodic-in-space solutions $(\vu,p)$ with period $L$, that is, 
	$\vu(t,{\bf x} + L{\bf e}_i) = \vu(t,{\bf x})$ and $p(t,{\bf x}+L{\bf e}_i)=p(t,{\bf x})$ and $\theta(t,{\bf x}+L{\bf e}_i)=\theta(t,{\bf x})$ 
	almost surely 
	and for any $(t, {\bf x})\in (0,T)\times \mathbb{R}^2$  and $i=1,2$, where 
	$\{\bf e_1,e_2\}$ denotes the canonical basis of $\mathbb{R}^2$.

	The Boussinesq system is widely used to model buoyancy-driven flows in various physical and geophysical applications. 
	It forms the foundation of Rayleigh-Bénard convection models in fluid layers \cite{chandrasekhar2013hydrodynamic, bodenschatz2000recent} 
	and is also applied to thermal convection in porous media \cite{nield2006convection}. 
	In geophysical fluid dynamics, the Boussinesq approximation is essential for describing large-scale atmospheric motions, 
	ocean circulation, and stratified flows, where density variations are small but dynamically significant \cite{hsia2007stratified, cushman2011introduction, vallis2017atmospheric}. 
	Because real-world systems often involve uncertainty, these equations have been studied not only in deterministic settings but also in stochastic frameworks \cite{duan2009large, chueshov2010stochastic, alonso2020well, law2014existence, lin2024global}. 
	In our setting, the Navier--Stokes equations are coupled with a temperature transport equation that includes diffusion, 
	and the system is driven by multiplicative noise. Here, $\vu$ denotes the fluid velocity, $\theta$ its temperature, and $p$ the pressure.
	
	The well-posedness of system \eqref{eq1.1} has been extensively studied. 
	In the presence of multiplicative noise, the existence of martingale solutions was established in \cite{law2014existence, duan2009large}. 
	A general abstract framework covering equations of the form \eqref{eq1.1} was developed in \cite{chueshov2010stochastic}, 
	where the existence and uniqueness of strong variational solutions were proved. 
	More recently, the case with transport-type noise was analyzed in \cite{alonso2020well, lin2024global}, 
	establishing the existence and uniqueness of strong solutions.
	
	Numerical methods for stochastic Stokes and Navier-Stokes equations have been extensively investigated over the past decades. 
	These include finite element methods \cite{brzezniak2013finite, FV2024, CP2012, Breit, BM2021, BM_space, feng2025full, vo2025high, FPL2021,li2025optimal}, 
	splitting schemes \cite{BBM14, FL2020, CHP2012}, and other approaches. 
	Among these, the mixed finite element method is particularly natural, as it approximates velocity and pressure simultaneously. 
	
	A major challenge in the numerical analysis of stochastic Navier-Stokes equations, 
	expected to be similar for stochastic Boussinesq equations, is the interaction between the multiplicative noise and the nonlinear term that prevents direct application of the classical Gronwall inequality used in deterministic analysis. 
	To address this, \cite{CHP2012} introduced a localization technique based on subsets of the sample space, allowing estimating error estimates for the mixed finite element method on these subsets. 
	However, this approach guarantees only \emph{convergence in probability}. 
	Full-moment error estimates were later obtained in \cite{bessaih2019strong, BM_space} by analyzing errors on the complements of these subsets, 
	leading to stronger convergence in certain $L^2$-type norms, though without polynomial rates in the case of multiplicative noise. 
	Although \cite{BM2021} established strong error estimates with polynomial rates, their analysis was restricted to constant additive noise. 
	Thus, strong convergence with polynomial order for multiplicative noise remained open. 
	Recently, \cite{feng2025full} derived full-moment strong error estimates for a time discretization of \eqref{eq1.1} using a stochastic Gronwall inequality, 
	obtaining polynomial rates of convergence under bounded multiplicative noise.
	
	Numerical methods for the stochastic Boussinesq equations have been considered for the first time in \cite{bessaih2022speed, bessaih2025rate}, although they only considered the fully implicit Euler method for time discretization of \eqref{eq1.1}. 
	Using a localization technique, error estimates for velocity and temperature in $L^2$ and $H^1$ norms were obtained, achieving convergence of order $O(k^{\alpha})$ for $\alpha \in (0,1/2)$. 
	However, no error estimates for the pressure were considered, and spatial discretization was not addressed, leaving a gap in the full numerical approximation of the stochastic Boussinesq equations.
	
	In this paper, we aim to fill this gap by studying a fully discrete mixed finite element scheme for system \eqref{eq1.1}. 
	To the best of our knowledge, this is the first study to combine both temporal and spatial discretizations within a unified framework. 
	Our approach provides a complete numerical approximation, capturing both the temporal evolution and spatial structure of the solutions. 
	
	Specifically, the fully discrete scheme in the Main Algorithm (Section~\ref{section 3}) employs the semi-implicit Euler--Maruyama method for time stepping, 
	effectively handling the stochastic terms, together with the standard mixed finite element method for spatial discretization, 
	which naturally approximates both velocity and pressure simultaneously. 
	Using the localization technique, we derive rigorous error estimates demonstrating convergence in probability. 
	The convergence of velocity and temperature is established in $L^2$ and $H^1$ norms, 
	and error estimates for pressure are derived in a time-averaged norm. 
	Overall, this work provides a comprehensive theoretical justification for the proposed scheme and confirms its effectiveness in approximating solutions of the stochastic Boussinesq equations.
	
	The remainder of the paper is organized as follows. 
	In Section~\ref{sec-2}, we introduce notations and preliminaries, including solution definitions and the well-posedness of \eqref{eq1.1}. 
	Section~\ref{section 3} presents the fully discrete mixed finite element time-stepping scheme, 
	derives stability estimates for velocity, pressure, and temperature, 
	and introduces the standard mixed finite element method for spatial discretization, with the MINI element as a prototypical example. 
	We also establish strong $L^2$ and $H^1$ error estimates in Theorem~\ref{Main_Theorem_error} using the localization technique and high-moment stability estimates. 
	Finally, Section~\ref{section 4} provides numerical experiments to validate the theoretical error estimates.

	\section{Preliminaries}\label{sec-2}
	\subsection{Notations}\label{sec-2.1}
	Standard function and space notation will be adopted in this paper. 
	We denote ${L}^p(D)$ and ${ H}^{k}(D)$ as the Lebesgue and Sobolev spaces of the functions that are periodic with period $L$ and have vanishing mean, while we also denote $\vL^p(D)$ and $\vH^{k}(D)$ as their vector-valued spaces. 
	$C$ denotes a generic constant that is independent of the mesh parameters $h$ and $k$.
	
	Let $(\Omega,\cF, \{\cF_t\},\mP)$ be a filtered probability space with the probability measure $\mP$, the 
	$\sigma$-algebra $\cF$ and the continuous  filtration $\{\cF_t\} \subset \cF$. For a random variable $v$ 
	defined on $(\Omega,\cF, \{\cF_t\},\mP)$,
	${\mathbb E}[v]$ denotes the expected value of $v$. 
	For a vector space $X$ with norm $\|\cdot\|_{X}$,  and $1 \leq p < \infty$, we define the Bochner space
	$\bigl(L^p(\Omega;X); \|v\|_{L^p(\Omega;X)} \bigr)$, where
	$\|v\|_{L^p(\Omega;X)}:=\bigl({\mathbb E} [ \Vert v \Vert_X^p]\bigr)^{\frac1p}$.
	We also define 
	\begin{align*}
		{\mathbb H} := \bigl\{{\bf v}\in  \vL^2(D) ;\,\div {\bf v}=0 \mbox{ in }D\, \bigr\}\, , \quad 
		{\mathbb V} :=\bigl\{{\bf v}\in  \vH^1(D) ;\,\div {\bf v}=0 \mbox{ in }D \bigr\}\, .
	\end{align*}
	
	We recall from \cite{Girault_Raviart86} the (orthogonal) Helmholtz projection 
	${\bf P}_{{\mathbb H}}: \vL^2(D) \rightarrow {\mathbb H}$  and define the Stokes operator ${\bf A} := -{\bf P}_{\mathbb H} \Delta: {\mathbb V} \cap \vH^2(D) \rightarrow {\mathbb H}$. 
	

	\subsection{Solution concepts and H\"older continuity estimates}\label{sub5.2.2}
	
	\smallskip
	
	First, we state clearly the assumptions on the diffusion functions $\vG_1$ and $G_2$. Suppose that 
	\smallskip
	
	\begin{enumerate}[(B1)]
		\item Let $\vG_1: \mH \rightarrow \mH$ be such that 
		\begin{align*}
			\|{\vG}_1(\vu)\|_{\vL^2} &\leq C_{\vG_1}(1 + \|\vu\|_{\vL^2})\qquad\forall \vu \in \mH,\\
			\|{\vG}_1(\vu) - \vG_1(\vv)\|_{\vL^2} &\leq C_{\vG_1}\|\vu-\vv\|_{\vL^2}\qquad\forall \vu,\vv \in \mH.
		\end{align*}
		
		\smallskip
		\item Let $\vG_1: \mV \rightarrow \mV$ be such that 
		\begin{align*}
			\|{\vG}_1(\vu)\|_{\vH^1} &\leq C_{\vG_1}(1 + \|\vu\|_{\vH^1})\qquad\forall \vu \in \mV,\\
			\|{\vG}_1(\vu) - \vG_1(\vv)\|_{\vH^1} &\leq C_{\vG_1}\|\vu-\vv\|_{\vH^1}\qquad\forall \vu,\vv \in \mV.
		\end{align*}
		
		\item Let $G_2: H^m(D) \rightarrow H^m(D)$ for $m =  0, 1$ be such that 
		\begin{align*}
			\|{G}_2(\theta)\|_{H^m} &\leq C_{G_2}(1 + \|\theta\|_{H^m})\qquad\forall \theta \in H^m(D),\\
			\|{G}_2(\theta) - G_2(\phi)\|_{H^m} &\leq C_{G_2}\|\theta-\phi\|_{H^m}\qquad\forall \theta, \phi \in H^m(D).
		\end{align*}

		\smallskip
		
		\item Let $\vG_1: \vH^2(D) \rightarrow \vH^2(D)$ be such that 
		\begin{align*}
			\|{\vG}_1(\vu)\|_{\vH^2} &\leq C_{\vG_1}(1 + \|\vu\|_{\vH^2})\qquad\forall \vu \in \vH^2(D).
		\end{align*}

		\smallskip
	\end{enumerate}
	
	\bigskip
	
	Next, we state the result about the existence and uniqueness of a weak pathwise solution of \eqref{eq1.1} in the following theorem. Its proof can be found in \cite{chueshov2010stochastic,duan2009large,bessaih2022speed}.
	
	\smallskip 
	
	\begin{theorem}\label{theorem2.2} 
		Given $(\Omega,\cF, \{\cF_t\},\mP)$, let $W_1$ and $W_2$ be ${\mathbb R}$-valued Wiener processes on it. 
		Suppose $(\vu_0,\theta_0)\in L^{2q}(\Omega; {\mathbb H}\times L^2(D))$ for $q \in \mathbb{N}$. Assume that $\vG_1$ and $G_2$ satisfy the conditions (B1)--(B3).	There exist a unique weak pathwise pair $\{(\vu(t), \theta(t)) ; 0\leq t\leq T\}$ of \eqref{eq1.1} such that ${\bf u} \in  L^2\bigl(\Omega; C(0,T; {\mathbb H}))\cap L^2\bigl(\Omega; L^2(0,T; {\mathbb V}))
		$ and $\theta \in L^2\bigl(\Ome; C(0,T; L^2(D))\bigr)\cap L^2\bigl(\Ome;L^2(0,T; H^1(D))$
		and satisfies $\mP$-a.s.~for all $t\in (0,T]$
		\begin{align}\label{equu2.8a}
			\bigl({\bf u}(t),  {\bf v} \bigr)& + \nu \int_0^t  \bigl(\nab {\bf u}(s), \nab {\bf v} \bigr) 
			\,  ds + \int_0^t \big([\vu(s)\cdot\nab]\vu(s),\vv\big)\, ds
			\\\nonumber
			&=({\bf u}_0, {\bf v}) + \int_{0}^t \bigl(\theta(s){\bf e}_2, \vv\bigr)\, ds + {  \Bigl(\int_0^t \vG_1(\vu(s))\, dW_1(s), {\bf v} \Bigr)}  \qquad\forall  \, {\bf v}\in {\mathbb V}, \\
			\bigl(\theta(t),  \phi \bigr)& + \mu \int_0^t  \bigl(\nab \theta(s), \nab \phi \bigr)\, ds 
			+ \int_0^t \big(\vu(s)\cdot\nab \theta(s),\phi\big)\, ds
			\\\nonumber
			&=(\theta_0, \phi) + {  \Bigl(\int_0^t G_2(\theta(s))\, dW_2(s), \phi \Bigr)}  \qquad \forall \phi \in H^1(D). 
		\end{align}
	\end{theorem}
	
	We also recall some properties of the convective term. 
	It is easy to check that 
	\begin{align*}
		\bigl([\vu\cdot\nab]\vv,\vv\bigr) &=0  \qquad\forall \vu\in \mV, \quad\forall\vv \in \vH^1(D).\\
		\bigl(\vu\cdot\nab\phi,\phi\bigr) &= 0 \qquad\forall \vu\in \mV, \quad\forall\phi\in H^1(D).
	\end{align*}

	Theorem \ref{theorem2.2} only gives the velocity $\mathbf{u}$ and the temperature $\theta$ for \eqref{eq1.1}, 
	its associated pressure $p$ is subtle to define. In \cite{LRS03}, the pressure solution is defined as a distribution, which is very difficult to analyze in error analysis later. So, we come up with a new definition of the pressure in the following theorem and refer the reader to \cite[Theorem 1.3]{FPL2021} for a similar proof.
	\begin{theorem}\label{thm2.2}
		Let $\{(\vu(t), \theta(t)) ; 0\leq t\leq T\}$ be the variational solution of \eqref{eq1.1}. There exists a unique adapted process 
		$P\in {L^2\bigl(\Omega; L^2(0,T; H^1(D)/\mathbb{R})\bigr)}$ such that $(\mathbf{u}, \theta, P)$ satisfies 
		$\mP$-a.s.~for all $t\in (0,T]$
		\begin{subequations}\label{New_variation}
			\begin{align}\label{New_variation_a}
				&\bigl({\bf u}(t),  {\bf v} \bigr) + \nu\int_0^t  \bigl(\nab {\bf u}(s), \nab {\bf v} \bigr) \, ds + \int_0^t \big([\vu(s)\cdot\nab]\vu(s),\vv\big)\, ds
				- \bigl(  \div \mathbf{v}, P(t) \bigr) \\
				&\qquad=({\bf u}_0, {\bf v}) + \int_{0}^t \bigl(\theta(s){\bf e}_2, \vv\bigr)\, ds
				+  {\int_0^t  \bigl( {\bf G}_1(\vu(s)), {\bf v} \bigr)\, dW_1(s)}  \qquad \forall  \, {\bf v}\in \vH^1(D)\, , \nonumber \\ \label{New_variation_b}
				&	\bigl(\theta(t),  \phi \bigr) + \mu \int_0^t  \bigl(\nab \theta(s), \nab \phi \bigr)\, ds 
				+ \int_0^t \big(\vu(s)\cdot\nab \theta(s),\phi\big)\, ds
				\\\nonumber
				&\qquad=(\theta_0, \phi) + {  \Bigl(\int_0^t G_2(\theta(s))\, dW_2(s), \phi \Bigr)}  \qquad\forall \phi \in H^1(D). \\
				&\bigl(\div {\bf u}, q \bigr) =0 \qquad\forall \, q\in L^2(D)/\mathbb{R}  .  \label{New_variation_c}
			\end{align}
		\end{subequations}
	\end{theorem}

	Next, we state higher regularity results for $(\vu, \theta)$, which are used to obtain our error estimates later.

	\begin{lemma}\label{stability_pdes} 
		Let $(\vu, \theta)$ be the solution from Theorem \ref{theorem2.2}. Then we have 
		\begin{enumerate}[{\rm (a)}]
			\item Assume that $\vu_0 \in L^{2q}\bigl(\Ome; \mV\bigr)$ and $\theta_0 \in L^{2q}(\Ome; L^2(D))$ for $q=1$ or $q\geq 2$. Additionally, suppose that $\vG_1$ and $G_2$ satisfy the conditions (B1)--(B3). Then there exists $C= C(T,q,\vu_0,\theta_0)>0$ such that
			\begin{align*}
				\mE\left[\sup_{t\in[0,T]}\|\nab\vu(t)\|^{2q}_{\vL^2} + \int_{0}^T \|\vA\vu(t)\|^2_{\vL^2}[1+ \|\vA^{\frac12} \vu(t)\|^{2(q-1)}_{\vL^2}]\, dt\right] \leq C.
			\end{align*}
			\item Assume that $\vu_0 \in L^{2}\bigl(\Ome; H^2(D)\bigr)\cap L^{10q}(\Ome; \mV)$ and $\theta_0 \in L^{2q}(\Ome; L^2(D))$ for $q\geq 1$. Additionally, suppose that $\vG_1$ and $G_2$ satisfy the conditions (B1)--(B4). Then there exists $C= C(T,q,\vu_0,\theta_0)>0$ such that
			\begin{align*}
				\mE\left[\sup_{t\in[0,T]}\|\Delta\vu(t)\|^{2q}_{\vL^2} \right] \leq C.
			\end{align*}
			\item Assume that $\vu_0 \in L^{9q}\bigl(\Ome; \mV\bigr)$ and $\theta_0 \in L^{9q}(\Ome; H^1(D))$ for $q=1$ or $q\geq 2$. Additionally, suppose that $\vG_1$ and $G_2$ satisfy the conditions (B1)--(B3). Then there exists $C= C(T,q,\vu_0,\theta_0)>0$ such that
			\begin{align*}
				\mE\left[\sup_{t\in[0,T]}\|\nab\theta(t)\|^{2q}_{L^2} + \int_{0}^T \|\Delta\theta(t)\|^2_{L^2}\|\nab \theta(t)\|^{2(q-1)}_{L^2}\, dt\right] \leq C.
			\end{align*}
		\end{enumerate}
	\end{lemma}
	\begin{proof}
		The proofs of (a) and (c) can be found in \cite[Propositions 1, 2]{bessaih2022speed}. Additionally, the proof of (b) follows the same lines with the proof of \cite[Lemma 2 (c)]{Breit} with a little modification on bounding $\theta$ in  $L^2\bigl(\Ome;L^2([0,T]); H^1(D)\bigr)$, which can be easily obtained by using Theorem \ref{theorem2.2}. So, we left this proof as an exercise for the interested reader.
	\end{proof}
	
	%
	%

	\subsection{H\"older continuity of the variational solution}
	We state the following high moment H\"older continuity estimates for $(\vu,\theta)$.  
	\begin{lemma}\label{lemma2.2} Let $(\vu,P,\theta)$ be solution from Theorem \ref{thm2.2} and let $q\in [2,\infty)$ and $\alpha \in \left(0,\frac12\right)$. Assume that $\vG_1$ and $G_2$ satisfy the conditions (B1)--(B3). Then, we have
		\begin{enumerate}[{\rm (a)}]
			\item If $\vu_0\in L^{4q}(\Ome; \mV)$ and $\theta_0 \in L^{2q}(\Ome; L^2(D))$. Then, there exists $C = C(q,\vu_0,\theta_0)>0$ such that
			\begin{align}
				&	\mE\left[\|\vu(t) - \vu(s)\|^{2q}_{\vL^2}\right] \leq C|t-s|^q\qquad\forall t, s \in [0,T],\\
				&	\mE\left[\left|\sum_{n=1}^{M}\int_{t_{n-1}}^{t_{n}} \bigl[\|\vu(s) - \vu(t_n)\|^2_{\vH^1} + \|\vu(s) - \vu(t_{n-1})\|^2_{\vH^1}\bigr]\right|^q\right]
				\leq Ck^{2\alpha q}.
			\end{align}
			\item If $\vu_0\in L^{9q}(\Ome; \mV)$ and $\theta_0 \in L^{9q}(\Ome; H^1(D))$. Then, there exists $C = C(q,\vu_0,\theta_0)>0$ such that
			\begin{align*}
				\mE\left[\|\theta(t) - \theta(s)\|^{2q}_{L^2}\right] \leq C|t-s|^q\qquad\forall t, s \in [0,T].
			\end{align*}
			\item If $\vu_0\in L^{17q}(\Ome; \mV)$ and $\theta_0 \in L^{17q}(\Ome; L^2(D))$. Then, there exists $C = C(q,\vu_0,\theta_0)>0$ such that
			\begin{align*}
				\mE\left[\left|\sum_{n=1}^{M}\int_{t_{n-1}}^{t_{n}} \bigl[\|\theta(s) - \theta(t_n)\|^2_{H^1} + \|\theta(s) - \theta(t_{n-1})\|^2_{H^1}\bigr]\right|^q\right] \leq Ck^{2\alpha q}.
			\end{align*}
			\item If $\vu_0 \in L^{2}\bigl(\Ome; H^2(D)\bigr)\cap L^{20}(\Ome; \mV)$ and $\theta_0 \in L^{2}(\Ome; L^2(D))$. Additionally, assume that $\vG_1$ satisfies the condition (B4). Then, there exists $C = C(q,\vu_0,\theta_0)>0$ such that
			\begin{align*}
				\mE\left[\|\nab(P(t) - P(s))\|^2_{L^2}\right] \leq C(t-s)^2.
			\end{align*}
		\end{enumerate}
	\end{lemma}
	\begin{proof}
		The proofs of (a), (b), (c) can be found in \cite[Propositions 4, 5 ]{bessaih2022speed}. Here we present the proof of (d) only. Using integrating by parts, we obtain the strong form of \eqref{New_variation_a}
		\begin{align}\label{eq2.9}
			\nab(P(t) - P(s)) &= - (\vu(t) - \vu(s))+ \nu \int_{s}^{t}\Delta \vu(\xi)\, d\xi - \int_{s}^{t}[\vu(\xi)\cdot\nab]\vu(\xi)\, d\xi \\\nonumber&\qquad+ \int_{s}^{t}\theta(s){\bf e}_2\, d\xi+ \int_{s}^t \vG_1(
			\vu(\xi))\, dW_1(\xi).
		\end{align}
		
		Testing \eqref{eq2.9} by $\nab(P(t) - P(s))$ and using the integration by parts together with \eqref{New_variation_c} and the divergence-free condition of $\vG_1$, we have
		\begin{align}\label{eq2.7}
			\|\nab(P(t) - P(s))\|^2_{\vL^2} &= \left(\int_{s}^{t}[\vu(\xi)\cdot\nab]\vu(\xi)\, d\xi, \nab(P(t) - P(s))\right) \\\nonumber&\qquad+ \left(\int_{s}^{t}\theta(s){\bf e}_2\, d\xi, \nabla(P(t) - P(s))\right).
		\end{align}
		
		Next, using the Cauchy-Schwarz inequality on \eqref{eq2.7}, we get
		\begin{align}\label{eq2.10}
			\|\nab(P(t) - P(s))\|^2_{\vL^2} &= 2\left\|\int_{s}^{t}[\vu(\xi)\cdot\nab]\vu(\xi)\, d\xi\right\|^2_{\vL^2} +2 \left\|\int_{s}^{t}\theta(s){\bf e}_2\, d\xi\right\|^2_{L^2}.
		\end{align}
		
		Taking the expectation to \eqref{eq2.10} and then using the Cauchy-Schwarz inequality, Lemma \ref{stability_pdes} Parts (a), (b) with $q =2$, and Theorem \ref{theorem2.2}, we obtain
		\begin{align}\label{eq2.8}
			\mE[\|\nab(P(t) - P(s))\|^2_{\vL^2}] &\leq  C|t-s|\int_{s}^t\mE[\|[\vu(\xi)\cdot\nab]\vu(\xi)\|^2_{\vL^2}]\, d\xi \\\nonumber
			&\qquad+ C|t-s|\int_{s}^{t}\mE[\|\theta(\xi)\|^2_{L^2}]\,d\xi\\\nonumber
			&\leq  C|t-s|\int_{s}^t\mE[\|\vA \vu(s)\|^2_{\vL^2}\|\nab \vu(\xi)\|^2_{\vL^2}]\, d\xi \\\nonumber
			&\qquad+ C|t-s|\int_{s}^{t}\mE[\|\theta(\xi)\|^2_{L^2}]\,d\xi\\\nonumber
			&\leq  C|t-s|^2\sup_{s\in [0,T]}\mE[\|\vA \vu(s)\|^2_{\vL^2}\|\nab \vu(\xi)\|^2_{\vL^2}] \\\nonumber
			&\qquad+ C|t-s|^2\sup_{s\in [0,T]}\mE[\|\theta(\xi)\|^2_{L^2}]\\\nonumber
			&\leq C(t-s)^2.
		\end{align}
		
	\end{proof}


	\bigskip
	
	\section{Mixed finite element method}\label{section 3} 
	\subsection{Formulation and stability of the fully discrete mixed finite element method}
	Let $\mathcal{T}_h$ be a quasi-uniform mesh of the domain $D \subset \mathbb{R}^2$ with mesh size $h > 0$. We consider the following finite element spaces:
	\begin{align*}
		\mH_h &= \bigl\{\vv_h \in {\bf C}(\overline{D}) \cap \vH^1(D);\,\vv_h \in [\mathcal{P}_{\ell}(K)]^2\qquad\forall K \in \mathcal{T}_h\bigr\},\\
		L_h &= \bigl\{\psi_h \in C(\overline{D})\cap \in L^2(D)/\mathbb{R};\, \psi_h \in \mathcal{P}_{m}(K)\qquad\forall K \in \mathcal{T}_h\bigr\},\\
		M_h &= \bigl\{\phi_h \in {C}(\overline{D}) \cap H^1(D);\,\phi_h \in \mathcal{P}_{j}(K)\qquad\forall K \in \mathcal{T}_h\bigr\},
	\end{align*}  
	where $\mathcal{P}_{\ell}$, $\mathcal{P}_{m}$, and $\mathcal{P}_{j}$ are the spaces of piece-wise polynomials on $K$ with degree at most $\ell, m , j\geq 1$. To ensure the stability of later schemes, we require that the pair $\mH_h$ and $L_h$ satisfy the Ladyzhenskaya-Babuska-Brezzi (LBB)  (or inf-sup condition), which is now quoted:  there exists $\beta_1>0$ such that
	\begin{align}\label{inf-sup_discrete}
		\sup_{\pphi_h \in \mH_h, \pphi_h \neq 0} \frac{\bigl(\div \pphi_h,\psi_h\bigr)}{\|\pphi_h\|_{\mH_h}} \geq \beta_1\|\psi_h\|_{L^2}\qquad\forall \psi_h\in L_h,
	\end{align}
	where the constant $\beta_1$ is independent of $h$ (and $k$). 
	
	There are many well-known pairs of finite elements that satisfy the condition \eqref{inf-sup_discrete}. For example, the Taylor-Hood elements, the MINI element, etc. For presentation purposes, in this paper, we only consider the MINI element for the pair $(\mH_h, L_h)$ in the error analysis, while the piecewise linear polynomial is chosen for $M_h$.
	
	Next, we define the space of a weakly divergent-free velocity field as follows:
	\begin{align*}
		\mV_h := \bigl\{\pphi_h \in \mH_h;\, \bigl(\div \pphi_h, q_h\bigr) =0\qquad\forall q_h \in L_h\bigr\}.
	\end{align*}
	In general, $\mV_h$ is not a subspace of $\mV$.
	Let $\cQ_h: \vL^2 \rightarrow \mV_h$ denote the $L^2$-orthogonal projection, which is defined by
	\begin{align}\label{Q_property}
		\bigl(\vv - \cQ_h\vv, \pphi_h\bigr) = 0\qquad\forall \pphi_h\in \mV_h.
	\end{align}
	It is well-known \cite{Girault_Raviart86} that $\cQ_h$ satisfies the following estimates:
	\begin{align}\label{Q_property_a}
		\|\vv - \cQ_h\vv\|_{\vL^2} + h\|\nab(\vv - \cQ_h\vv)\|_{\vL^2} &\leq C h^2\|\vA\vv\|_{\vL^2}\qquad\forall \vv\in \mV\cap\vH^2(D),\\
		\label{Q_property_b}	\|\vv - \cQ_h\vv\|_{\vL^2} &\leq Ch\|\nab\vv\|_{\vL^2}\qquad\forall \vv\in \mV\cap\vH^1(D).
	\end{align} 
	Similarly, let $P_h: L^2_{per}(D) \rightarrow L_h$ and $B_h: L^2_{per}(D) \rightarrow M_h$ be the scalar $L^2$-orthogonal projections, defined by
	\begin{align}
		\label{P_property}	\bigl(\psi - P_h\psi, q_h\bigr) &= 0\qquad\forall q_h \in L_h,\\
		\label{B_property}	\bigl(\psi - B_h\psi, \phi_h\bigr) &= 0\qquad\forall \phi_h \in M_h,
	\end{align}
	then there also holds
	\begin{align}\label{P_property_a}
		\|\psi - P_h\psi\|_{L^2} &\leq Ch\|\nab \psi\|_{\vL^2}\qquad\forall \psi \in L^2(D)\cap H^1(D).\\
		\label{B_property_b}	\|v - B_hv\|_{L^2} + h\|\nab(v - B_h v)\|_{L^2} &\leq C h^2\|v\|_{H^2}\qquad\forall v\in H^2(D).
	\end{align}
	
	Next, we also introduce convenient trilinear forms \cite{Temam}
	\begin{align}
		{b}(\vu,\vv,\vw) &= \bigl([\vu\cdot\nab]\vv,\vw\bigr) + \frac12 \bigl([\div \vu]\vv,\vw\bigr)\qquad\forall \vu, \vv, \vw \in \vH^1(D).\\
		\tilde{b}(\vu,v,w) &= \bigl(\vu\cdot\nab v, w\bigr) + \frac12 \bigl([\div \vu] v, w\bigr)\qquad\forall \vu \in \vH^1(D), \,\forall v, w \in H^1(D).
	\end{align}
	which is anti-symmetric in the sense that
	\begin{align}
		{b}(\vu,\vv,\vw) &= - {b}(\vu,\vw,\vv) \qquad\forall\vu,\vv,\vw \in \vH^1(D).\\
		\tilde{b}(\vu,v,w) &= - \tilde{b}(\vu,w,v) \qquad\forall\vu\in \vH^1(D),\,\forall v,w \in H^1(D).
	\end{align}
	Therefore,\begin{align}\label{b_property}
		{b}(\vu, \pphi,\pphi) &= 0\qquad\forall\vu,\pphi \in \vH^1(D).\\
		\label{bb_property}		\tilde{b}(\vu, \phi,\phi) &= 0\qquad\forall\vu \in \vH^1(D), \, \forall \phi \in H^1(D).
	\end{align}

	Our fully discrete mixed finite element method is defined by the following algorithm. 
	
	\medskip
	
	\textbf{Main Algorithm.} Let $(\vu_h^0, \theta_h^0) \in \vH_h \times M_h$. Find $\displaystyle \bigl(\vu_h^{n+1},p_h^{n+1}, \theta_h^{n+1}\bigr) \in L^2(\Ome;\mH_h\times L_h\times M_h)$ such that there holds $\mP$-a.s.
	\begin{subequations}\label{Main_Scheme}
		\begin{align}
			\label{Main_Scheme_a}	&\bigl(\vu_h^{n+1} - \vu_h^n,\pphi_h\bigr) + \nu k \bigl(\nab \vu_h^{n+1},\nab\pphi_h\bigr) + k\, b(\vu_h^{n}, \vu_h^{n+1},\pphi_h)- k\bigl(p_h^{n+1},\div \pphi_h\bigr) \\\nonumber&\qquad= k \bigl(\theta_h^n{\bf e}_2, \pphi_h\bigr)+\bigl(\vG_1(\vu_h^n)\Delta W_{1,n+1},\pphi_h\bigr)\quad\forall \pphi_h \in \vH_h,\\
			\label{Main_Scheme_b}	&\bigl(\theta_h^{n+1} - \theta_h^n,\varphi_h\bigr) + \mu k \bigl(\nab \theta_h^{n+1},\nab\varphi_h\bigr) + k\, \tilde{b}(\vu_h^{n+1}, \theta_h^{n+1},\varphi_h)  \\\nonumber
			&\qquad= \bigl(G_2(\theta_h^n)\Delta W_{2,n+1},\varphi_h\bigr) \quad\forall \varphi_h\in M_h,\\
			&(\div\vu_h^{n+1},\psi_h\bigr) =0 \quad \forall \psi_h \in L_h,\label{Main_Scheme_c}
		\end{align}
	\end{subequations}
	where $\Delta W_{i, n+1} = W(t_{n+1}) - W(t_n) \sim \mathcal{N}(0,k)$ for $i = 1,2$.
	
	\medskip

	Next, we derive the following stability estimates for $(\vu_h^n, \theta_h^n)$.
	
	\begin{lemma}\label{stability_FEMs}
		Let $1 \leq q \leq 3$ and $(\vu_h^0,\theta^0_h) \in L^{2^q}(\Ome;\mH_h\times M_h)$. Let $\{(\vu_h^n, \theta_h^n)\}$ be the approximate solution from the scheme \eqref{Main_Scheme}. Then, there hold
		\begin{enumerate}[{\rm (a)}]
			\item $\displaystyle \mE\biggl[\max_{1\leq n \leq M}\|\vu_h^n\|^{2^q}_{\vL^2} + \nu k\sum_{n=1}^M \|\vu_h^n\|^{2^{q}-2}_{\vL^2}\|\nab\vu_h^n\|^2_{\vL^2}\biggr] \leq C_1,$
			\item $\displaystyle \mE\biggl[\max_{1\leq n \leq M}\|\theta_h^n\|^{2^q}_{L^2} + \mu k\sum_{n=1}^M \|\theta_h^n\|^{2^{q}-2}_{L^2}\|\nab\theta_h^n\|^2_{L^2}\biggr] \leq C_2,$
		\end{enumerate}
		where $C_1 = C(T,q,\vu_h^0)$ and $C_2 = C(T,q,\theta_h^0)$.
	\end{lemma}
	\begin{proof}
		First, the proof of (a) can be found in \cite[Lemma 3.3]{brzezniak2013finite} with the help of using (b) to control the terms related to $\theta_h^n$. 
		
		We only need to prove (b). Taking $\varphi_h = \theta_h^{n+1}$ in \eqref{Main_Scheme_b} and using the identity $2a(a-b) = a^2 - b^2 + (a-b)^2$, and also using the property \eqref{bb_property} we obtain
		\begin{align*}
			&	\frac12\left[\|\theta_h^{n+1}\|^2_{L^2} - \|\theta_h^n\|^2_{L^2} + \|\theta_h^{n+1} - \theta_h^n\|^2_{L^2}\right] + \mu\|\nab \theta_h^{n+1}\|^2_{L^2} \\\nonumber&= \bigl(G_2(\theta_h^n)\Delta W_{2,n+1},\theta_h^{n+1} - \theta_h^n\bigr) +  \bigl(G_2(\theta_h^n)\Delta W_{2,n+1},\theta_h^n\bigr)\\\nonumber
			&\leq \frac{1}{4}\|\theta_h^{n+1} - \theta_h^n\|^2_{L^2} + \|G_2(\theta_h^n)\Delta W_{2,n+1}\|^2_{L^2} +  \bigl(G_2(\theta_h^n)\Delta W_{2,n+1},\theta_h^n\bigr),
		\end{align*}
		
		which implies that
		\begin{align}\label{eq3.16}
			&	\frac12\left[\|\theta_h^{n+1}\|^2_{L^2} - \|\theta_h^n\|^2_{L^2} + \frac12\|\theta_h^{n+1} - \theta_h^n\|^2_{L^2}\right] + \mu k\|\nab \theta_h^{n+1}\|^2_{L^2} \\\nonumber
			&\leq  \|G_2(\theta_h^n)\Delta W_{2,n+1}\|^2_{L^2} +  \bigl(G_2(\theta_h^n)\Delta W_{2,n+1},\theta_h^n\bigr).
		\end{align}
		
		Taking the expectation and using the fact that $\mE\left[\bigl(G_2(\theta_h^n)\Delta W_{n+1},\theta_h^n\bigr)\right] = 0$, and then applying the summation $\sum_{n = 0}^{\ell}$ for any $0\leq \ell \leq M-1$ to \eqref{eq3.16}, we obtain
		\begin{align}\label{eq3.17}
			&	\frac12\mE\left[\|\theta_h^{\ell+1}\|^2_{L^2}\right] + \frac14\sum_{n = 0}^{\ell}\mE\left[\|\theta_h^{n+1} - \theta_h^n\|^2_{L^2}\right] + \mu k \sum_{n = 0}^{\ell} \mE\left[\|\nab \theta_h^{n+1}\|^2_{L^2}\right]\\\nonumber
			&\leq \frac12\mE\left[\|\theta_h^0\|^2_{L^2}\right] + C_{G_2}k\sum_{n=0}^{\ell} \mE\left[\|\theta_h^n\|^2_{L^2}\right],
		\end{align}
		where the second term on the right-hand side of \eqref{eq3.17} is obtained by using (B2) and the independence property of $\Delta W_{n+1}$. 
		
		Next, applying the discrete Gronwall inequality to \eqref{eq3.17} and then taking the maximum, we obtain
		\begin{align}\label{eq3.18}
			&	\max_{1\leq \ell \leq M}\mE\left[\|\theta_h^{\ell}\|^2_{L^2}\right] + \sum_{n = 0}^{M-1}\mE\left[\|\theta_h^{n+1} - \theta_h^n\|^2_{L^2}\right] + \mu k \sum_{n = 0}^{M} \mE\left[\|\nab \theta_h^{n}\|^2_{L^2}\right]\\\nonumber
			&\leq 2\mE\left[\|\theta_h^0\|^2_{L^2}\right]\exp(4C_{G_2}T).
		\end{align}
		
		Now, we use \eqref{eq3.18} as a tool \cite{brzezniak2013finite,vo2025high} to prove (b) in the case $q =1$. To do that, first applying the summation $\sum_{n = 0}^{\ell}$to \eqref{eq3.16} and $\max_{1\leq \ell \leq M-1}$ to \eqref{eq3.16}, and then taking the expectation, we obtain
		\begin{align}\label{eq3.19}
			&	\mE\left[\max_{1\leq \ell \leq M}\|\theta_h^{\ell}\|^2_{L^2}\right] + \sum_{n = 0}^{M}\mE\left[\|\theta_h^{n+1} - \theta_h^n\|^2_{L^2}\right] + \mu k \sum_{n = 0}^{M} \mE\left[\|\nab \theta_h^{n}\|^2_{L^2}\right]\\\nonumber
			&\leq  2\mE\left[\|\theta_h^0\|^2_{L^2}\right]+4\sum_{n = 0}^{M-1}\mE\left[\|G_2(\theta_h^n)\Delta W_{2,n+1}\|^2_{L^2}\right] \\\nonumber
			&\qquad+ 4\mE\left[ \max_{1\leq \ell \leq M-1}\sum_{n = 0}^{\ell}\bigl(G_2(\theta_h^n)\Delta W_{2,n+1},\theta_h^n\bigr)\right]\\\nonumber
			&\leq  2\mE\left[\|\theta_h^0\|^2_{L^2}\right]+4C_{G_2}k\sum_{n = 0}^{M-1}\mE\left[\|\theta_h^n\|^2_{L^2}\right] + 4\mE\left[ \max_{1\leq \ell \leq M-1}\sum_{n = 0}^{\ell}\bigl(G_2(\theta_h^n)\Delta W_{2,n+1},\theta_h^n\bigr)\right],
		\end{align}
		
		where the second term on the last inequality in the right-hand side of \eqref{eq3.19} is obtained by using (B2) and the independent property of $\Delta W_{n+1}$. Next, we can use \eqref{eq3.18} to bound the second term, while using the Burkholder-Davis-Gundy inequality, we estimate the last term on the right-hand side of \eqref{eq3.19} as follows:
		\begin{align*}
			&\mE\left[ \max_{1\leq \ell \leq M-1}\sum_{n = 0}^{\ell}\bigl(G_2(\theta_h^n)\Delta W_{2,n+1},\theta_h^n\bigr)\right] \leq \mE\left[\left(k\sum_{n = 0}^M\|G_2(\theta_h^n)\|^2_{L^2}\|\theta_h^{n}\|^2_{L^2}\right)^{\frac12}\right]\\\nonumber
			&\leq \frac14\mE\left[\max_{1\leq n \leq M}\|\theta_h^n\|^2_{L^2}\right] + C_{G_2}k\sum_{n = 0}^M \mE\left[\|\theta_h^n\|^2_{L^2}\right].
		\end{align*}
		
		With this and \eqref{eq3.19}, and \eqref{eq3.18}, we obtain (b) for $q=1$. 
		
		Next, we prove (b) with $q=2$. To do that, we multiply \eqref{eq3.16} with $\|\theta_h^{n+1}\|^2_{L^2}$ and also use the identity $2a(a-b) = a^2 - b^2 + (a-b)^2$
		to obtain
		\begin{align}
			&	\frac14\left[\|\theta_h^{n+1}\|^4_{L^2} - \|\theta_h^n\|^4_{L^2}\right] +\frac14\left(\|\theta_h^{n+1}\|^2_{L^2} - \|\theta_h^n\|^2_{L^2}\right)^2 + \frac14\|\theta_h^{n+1} - \theta_h^n\|^2_{L^2}\|\theta^{n+1}_h\|^2_{L^2}\\\nonumber
			&\qquad + \mu k\|\nab \theta_h^{n+1}\|^2_{L^2}\|\theta_h^{n+1}\|^2_{L^2} \\\nonumber
			&\leq  \|G_2(\theta_h^n)\Delta W_{2,n+1}\|^2_{L^2}\|\theta_h^{n+1}\|^2_{L^2} +  \bigl(G_2(\theta_h^n)\Delta W_{2,n+1},\theta_h^n\bigr)\|\theta_h^{n+1}\|^2_{L^2}\\\nonumber
			&=  \|G_2(\theta_h^n)\Delta W_{2,n+1}\|^2_{L^2}(\|\theta_h^{n+1}\|^2_{L^2} - \|\theta_h^n\|^2_{L^2}) + \|G_2(\theta_h^n)\Delta W_{2,n+1}\|^2_{L^2}\|\theta_h^n\|^2_{L^2} \\\nonumber
			&\qquad+  \bigl(G_2(\theta_h^n)\Delta W_{2,n+1},\theta_h^n\bigr)(\|\theta_h^{n+1}\|^2_{L^2} - \|\theta_h^n\|^2_{L^2}) +  \bigl(G_2(\theta_h^n)\Delta W_{2,n+1},\theta_h^n\bigr) \|\theta_h^n\|^2_{L^2}\\\nonumber
			&\leq \frac{1}{8}\left(\|\theta_h^{n+1}\|^2_{L^2} - \|\theta_h^n\|^2_{L^2}\right)^2_{L^2} +  4\|G_2(\theta_h^n)\Delta W_{2,n+1}\|^4_{L^2} + \|G_2(\theta_h^n)\Delta W_{2,n+1}\|^2_{L^2}\|\theta_h^n\|^2_{L^2} \\\nonumber
			&\qquad+  4\|G_2(\theta_h^n)\Delta W_{2,n+1}\|^2_{L^2}\|\theta_h^n\|^2_{L^2} +  \bigl(G_2(\theta_h^n)\Delta W_{2,n+1},\theta_h^n\bigr) \|\theta_h^n\|^2_{L^2},
		\end{align}
		which implies that
		\begin{align}\label{eq3.21}
			&	\frac14\left[\|\theta_h^{n+1}\|^4_{L^2} - \|\theta_h^n\|^4_{L^2}\right] +\frac18\left(\|\theta_h^{n+1}\|^2_{L^2} - \|\theta_h^n\|^2_{L^2}\right)^2 + \frac14\|\theta_h^{n+1} - \theta_h^n\|^2_{L^2}\|\theta^{n+1}_h\|^2_{L^2}\\\nonumber
			&\qquad + \mu k\|\nab \theta_h^{n+1}\|^2_{L^2}\|\theta_h^{n+1}\|^2_{L^2} \\\nonumber
			&\leq   4\|G_2(\theta_h^n)\Delta W_{2,n+1}\|^4_{L^2} + 5\|G_2(\theta_h^n)\Delta W_{2,n+1}\|^2_{L^2}\|\theta_h^n\|^2_{L^2}  \\\nonumber&\qquad+  \bigl(G_2(\theta_h^n)\Delta W_{2,n+1},\theta_h^n\bigr) \|\theta_h^n\|^2_{L^2}.
		\end{align}
		
		Taking the expectation and then applying the summation $\sum_{n = 0}^{\ell}$ for any $0\leq \ell \leq M-1$ to \eqref{eq3.21} we obtain
		\begin{align}\label{eq.3.22}
			&	\frac14\mE\left[\|\theta_h^{\ell+1}\|^4_{L^2}\right] +\frac18\sum_{n = 0}^{\ell}\mE\left[\left(\|\theta_h^{n+1}\|^2_{L^2} - \|\theta_h^n\|^2_{L^2}\right)^2\right] + \mu k\sum_{n = 0}^{\ell}\mE\left[\|\nab \theta_h^{n+1}\|^2_{L^2}\|\theta_h^{n+1}\|^2_{L^2}\right] \\\nonumber
			&\leq   	\frac14\mE\left[\|\theta_h^{0}\|^4_{L^2}\right] +9C_{G_2}\sum_{n = 0}^{\ell}\mE\left[\theta_h^n\|^4_{L^2}\right] + 0
		\end{align}
		
		Next, applying the discrete Gronwall inequality to \eqref{eq.3.22} and then taking the maximum, we obtain
		\begin{align}\label{eq.3.23}
			&	\max_{1\leq \ell \leq M}\mE\left[\|\theta_h^{\ell}\|^4_{L^2}\right] + \sum_{n = 0}^{M-1}\mE\left[\left(\|\theta_h^{n+1}\|^2_{L^2} - \|\theta_h^n\|^2_{L^2}\right)^2\right] + \mu k\sum_{n = 0}^{M}\mE\left[\|\nab \theta_h^{n}\|^2_{L^2}\|\theta_h^{n}\|^2_{L^2}\right] \\\nonumber
			&\leq   	2\mE\left[\|\theta_h^{0}\|^4_{L^2}\right]\exp(72C_{G_2}T).
		\end{align}
		
		Then, we can obtain the inequality in (b) with $q=2$ by using \eqref{eq.3.23} as a tool. This step is similar to the process of obtaining the inequality with $q=1$. So, we omit the details to save space. 
		
		Similarly, we can obtain the inequality in (b) with $q =3$ by multiplying \eqref{eq3.21} with $\|\theta_h^{n+1}\|^4_{L^2}$ and proceeding similarly to the cases $q = 1, 2$. So, we leave the details as an exercise for the interested reader.
		
		The proof is complete.
		
	\end{proof}

	\subsection{Error estimates}
	In this part, we state and prove the error estimates of the fully discrete mixed finite element approximations. To control the nonlinearity, we introduce the following sequence of subsets of the sample space
	\begin{align}
		{\Omega}_{\rho, m} := \left\{\omega \in \Omega; \, \sup_{t \leq t_m} \left(\|\vu(t)\|^2_{\vH^1} + \|\vu(t)\|^4_{\vH^1} +\|\theta(t)\|^2_{H^1} + \|\theta(t)\|^4_{H^1}\right) \leq \rho \right\},
	\end{align}
	where $(u,\theta)$ is the variational solution from Theorem \ref{theorem2.2} and for some $\rho>0$ specified later. We observe that ${\Omega}_{\rho,0} \supset {\Omega}_{\rho,1} \supset ... \supset {\Omega}_{\rho,\ell}$. 
	
	\begin{remark}
		In the error estimate stated in Theorem~\ref{Main_Theorem_error}, we choose
		\begin{align}
			\rho(k) := \frac{\ln(\ln(1/k))}{\tilde{C}T},
		\end{align}
		where
		\[
		\tilde{C} := \frac{4C_e^2}{\nu} + \frac{1024 C_e^4}{\nu^2},
		\]
		and $C_e$ denotes the constant appearing in the Ladyzhenskaya inequality in two dimensions.
		
		With this choice of $\rho(k)$, an application of the Markov inequality together with Lemma~\ref{stability_pdes}(a) and (c) yields
		\begin{align*}
			\mP(\Omega_{\rho,M}^c)
			&\le \frac{\tilde{C}T}{\ln(\ln(1/k))}
			\mE\!\left[
			\sup_{0\le t\le T}
			\Big(
			\|\vu(t)\|_{\vH^1}^2 + \|\vu(t)\|_{\vH^1}^4
			+ \|\theta(t)\|_{H^1}^2 + \|\theta(t)\|_{H^1}^4
			\Big)
			\right].
		\end{align*}
		Since the expectation on the right-hand side is finite, we conclude that
		\[
		\mP(\Omega_{\rho,M}^c) \longrightarrow 0
		\qquad \text{as } k \to 0,
		\]
		and hence
		\[
		\mP(\Omega_{\rho,M}) \longrightarrow 1
		\qquad \text{as } k \to 0.
		\]
		
		Therefore, the convergence of the numerical solutions implied by Theorem~\ref{Main_Theorem_error} is convergence in probability.
	\end{remark}

	\begin{theorem}\label{Main_Theorem_error} 
		Let $(u,\theta)$ be the variational solution to \eqref{New_variation} and $\{(u_h^{n},\theta_h^n)\}_{n=1}^M$ be generated by \eqref{Main_Scheme}. Let $u_0 \in L^2(\Omega; H^2(D))\cap L^{34}(\Ome;\mV)$ and $\theta_0\in L^{9}(\Omega; H^1(D))\cap L^{34}(\Omega; L^2(D))$. Assume that $\vG_1$ and $G_2$ satisfy the conditions (B1)--(B4). Let $\alpha \in \left(0,\frac{1}{2}\right)$. Then, there holds
		\begin{align}\label{eq3.30}
			&\left(\max_{1\leq n\leq M}\mE\left[\mathbf{1}_{{\Omega}_{\rho, M}}\|\vu(t_n) - \vu_h^n\|^{2}_{\vL^2}\right]\right)^{\frac{1}{2}} + \left(\max_{1\leq n\leq M}\mE\left[\mathbf{1}_{{\Omega}_{\rho, M}}\|\theta(t_n) - \theta_h^n\|^{2}_{L^2}\right]\right)^{\frac{1}{2}} \\\nonumber
			&\qquad+ \left(\mE\left[\mathbf{1}_{\Omega_{\rho, M}} k\sum_{n = 1}^M\left(\nu\|\nab(\vu(t_n) - \vu_h^n)\|^2_{\vL^2} + \mu\|\nab(\theta(t_n) - \theta_h^n)\|^2_{L^2}\right)\right]\right)^{\frac{1}{2}} \\\nonumber
			&\leq C\sqrt{\ln(1/k)}\left(k^{\alpha} + h\right),
		\end{align}
		where $C= C(u_0,\theta_0, T)$ is a positive constant. 
	\end{theorem}

	\begin{proof}
		First, denote $\ve^n := \vu(t_n) - \vu_h^n$ and $e_{\theta}^n := \theta(t_n) - \theta_h^n$	for any $0\leq n \leq M$, and $\tilde{p}(t) = \frac{P(t) - P(t-k)}{k}$ for any $t \in [0,T]$. Subtracting \eqref{Main_Scheme_a} to \eqref{New_variation_a}, we get the following error equations:
		\begin{align}\label{eq_error_u}
			&	\left(\ve^{n+1} - \ve^n, \pphi_h\right) + \nu k\left(\nab\ve^{n+1}, \nab \pphi_h\right) \\\nonumber
			&= \nu \int_{t_n}^{t_{n+1}} \left(\nab(\vu(t_{n+1}) - \vu(s)),\nab \pphi_h\right)\, ds + k\left(\tilde{p}(t_{n+1}) - p_h^{n+1}, \div \pphi_h\right)\\\nonumber
			&\qquad -\left[\int_{t_n}^{t_{n+1}}{b}(\vu(s),\vu(s),\pphi_h)\, ds - k{b}(\vu_h^{n}, \vu_h^{n+1},\pphi_h)  \right] \\\nonumber
			&\qquad + k \left(e_{\theta}^n{\bf e}_2,\pphi_h\right) + \int_{t_{n}}^{t_{n+1}}\left((\theta(s) - \theta(t_n)){\bf e}_2,\pphi_h\right)\, ds \\\nonumber
			&\qquad+ \int_{t_n}^{t_{n+1}}\left(\vG_1(\vu(s)) , \pphi_h\right)\, dW_1(s) - \left(\vG_1(\vu_h^n)\Delta W_{1,n+1}, \pphi_h\right)\\\nonumber
			&= \nu \int_{t_n}^{t_{n+1}} \left(\nab(\vu(t_{n+1}) - \vu(s)),\nab \pphi_h\right)\, ds + k\left(\tilde{p}(t_{n+1}) - p_h^{n+1}, \div \pphi_h\right)\\\nonumber
			&\qquad +\int_{t_n}^{t_{n+1}}\left[{b}(\vu(t_{n}),\vu(t_{n+1}),\pphi_h) - {b}(\vu(s), \vu(s),\pphi_h) \right] \, ds
			\\\nonumber
			&\qquad -k\left[{b}(\vu(t_{n}),\vu(t_{n+1}),\pphi_h) - {b}(\vu_h^{n}, \vu_h^{n+1},\pphi_h)  \right]\\\nonumber
			&\qquad + k \left(e_{\theta}^n{\bf e}_2,\pphi_h\right) + \int_{t_{n}}^{t_{n+1}}\left((\theta(s) - \theta(t_n)){\bf e}_2,\pphi_h\right)\, ds \\\nonumber
			&\qquad+ \left(\int_{t_n}^{t_{n+1}}\left(\vG_1(\vu(s)) - \vG_1(\vu(t_n))\right)\, dW_1(s), \pphi_h\right)\\\nonumber
			&\qquad+\left(\left(\vG_1(\vu(t_n))-\vG_1(\vu_h^n)\right)\Delta W_{1,n+1},\pphi_h\right).
		\end{align}
		
		Similarly, subtracting \eqref{New_variation_b} to \eqref{Main_Scheme_b}, we also get
		\begin{align}\label{eq_error_theta}
			&	\left(e_{\theta}^{n+1} - e_{\theta}^n, \varphi_h\right) + \mu k\left(\nab e_{\theta}^{n+1}, \nab \varphi_h\right) \\\nonumber
			&= \mu \int_{t_n}^{t_{n+1}} \left(\nab(\theta(t_{n+1}) - \theta(s)),\nab \varphi_h\right)\, ds \\\nonumber
			&\qquad -\left[\int_{t_n}^{t_{n+1}}\tilde{b}(\vu(s),\theta(s),\varphi_h)\, ds - k\tilde{b}(\vu_h^{n+1}, \theta_h^{n+1},\varphi_h)  \right] \\\nonumber
			&\qquad+ \int_{t_n}^{t_{n+1}}\left(G_2(\theta(s)) , \varphi_h\right)\, dW_2(s) - \left(G_2(\theta_h^n)\Delta W_{2,n+1}, \varphi_h\right)\\\nonumber
			&= \mu \int_{t_n}^{t_{n+1}} \left(\nab(\theta(t_{n+1}) - \theta(s)),\nab \varphi_h\right)\, ds \\\nonumber
			&\qquad + \int_{t_n}^{t_{n+1}}\left[\tilde{b}(\vu(t_n),\theta(t_{n+1}), \varphi_h)- \tilde{b}(\vu(s),\theta(s),\varphi_h)\right]\, ds    \\\nonumber
			&\qquad -k\left[\tilde{b}(\vu(t_{n+1}),\theta(t_{n+1}),\varphi_h)\, ds - \tilde{b}(\vu_h^{n+1}, \theta_h^{n+1},\varphi_h)  \right] \\\nonumber
			&\qquad+ \int_{t_n}^{t_{n+1}}\left(G_2(\theta(s)) , \varphi_h\right)\, dW_2(s) - \left(G_2(\theta_h^n)\Delta W_{2,n+1}, \varphi_h\right)\\\nonumber
			&\qquad+\left(\left(G_2(\theta(t_n))- G_2(\theta_h^n)\right)\Delta W_{2,n+1},\varphi_h\right).
		\end{align}

		Adding \eqref{eq_error_u} to \eqref{eq_error_theta}, and then taking $\pphi_h = \cQ_h\ve^{n+1}\in \mV_h$ and $\varphi_h = B_h e^{n+1}_{\theta} \in M_h$, and using \eqref{Q_property} and \eqref{B_property}, and $\left(p_h^{n+1},\div\cQ_h \ve^{n+1}\right) = 0$, we obtain
		\begin{align}\label{eq_3.19}
			&	\left(\cQ_h\ve^{n+1} - \cQ_h\ve^n, \cQ_h\ve^{n+1}\right) + \nu k\|\nab\ve^{n+1}\|^2_{\vL^2} \\\nonumber
			&\qquad\qquad+ \left(B_he_{\theta}^{n+1} - B_he_{\theta}^n, B_he_{\theta}^{n+1}\right) + \mu k\|\nab e_{\theta}^{n+1}\|^2_{L^2}\\\nonumber
			&=\nu k\left(\nab \ve^{n+1},\nab\left(\vu(t_{n+1})-  \cQ_h\vu(t_{n+1})\right)\right) \\\nonumber
			&\qquad+ \nu \int_{t_n}^{t_{n+1}} \left(\nab(\vu(t_{n+1}) - \vu(s)),\nab \cQ_h\ve^{n+1}\right)\, ds + k\left(\tilde{p}(t_{n+1}), \div \cQ_h\ve^{n+1}\right)\\\nonumber
			&\qquad +\int_{t_n}^{t_{n+1}}\left[{b}(\vu(t_{n}),\vu(t_{n+1}),\cQ_h\ve^{n+1}) - {b}(\vu(s), \vu(s),\cQ_h\ve^{n+1}) \right] \, ds
			\\\nonumber
			&\qquad -k\left[{b}(\vu(t_{n}),\vu(t_{n+1}),\cQ_h\ve^{n+1}) - {b}(\vu_h^{n}, \vu_h^{n+1}, \cQ_h\ve^{n+1})  \right]\\\nonumber
			&\qquad + k \left(e_{\theta}^n{\bf e}_2, \cQ_h\ve^{n+1}\right) + \int_{t_{n}}^{t_{n+1}}\left((\theta(s) - \theta(t_n)){\bf e}_2, \cQ_h\ve^{n+1}\right)\, ds \\\nonumber
			&\qquad+ \left(\int_{t_n}^{t_{n+1}}\left(\vG_1(\vu(s)) - \vG_1(\vu(t_n))\right)\, dW_1(s), \cQ_h\ve^{n+1}\right)\\\nonumber
			&\qquad+\left(\left(\vG_1(\vu(t_n))-\vG_1(\vu_h^n)\right)\Delta W_{1,n+1},\cQ_h\ve^{n+1}\right)\\\nonumber
			&\qquad+\mu k\left(\nab e_{\theta}^{n+1},\nab\left(\theta(t_{n+1})-  B_h\theta(t_{n+1})\right)\right) +\mu \int_{t_n}^{t_{n+1}} \left(\nab(\theta(t_{n+1}) - \theta(s)),\nab \varphi_h\right)\, ds \\\nonumber
			&\qquad + \int_{t_n}^{t_{n+1}}\left[\tilde{b}(\vu(t_{n+1}),\theta(t_{n+1}), B_h e^{n+1}_{\theta})- \tilde{b}(\vu(s),\theta(s),B_h e^{n+1}_{\theta})\right]\, ds    \\\nonumber
			&\qquad -k\left[\tilde{b}(\vu(t_{n+1}),\theta(t_{n+1}),B_h e^{n+1}_{\theta})\, ds - \tilde{b}(\vu_h^{n+1}, \theta_h^{n+1},B_h e^{n+1}_{\theta})  \right] \\\nonumber
			&\qquad+\left(\int_{t_n}^{t_{n+}} \left(G_2(\theta(s)) - G_2(\theta(t_n))\right)\, dW_2(s), B_h e^{n+1}_{\theta}\right)\\\nonumber
			&\qquad+\left(\left(G_2(\theta(t_n))- G_2(\theta_h^n)\right)\Delta W_{2,n+1},B_h e^{n+1}_{\theta}\right).
		\end{align}
		which can be rearranged as follows
		\begin{align}\label{eq_error}
			&	\left(\cQ_h\ve^{n+1} - \cQ_h\ve^n, \cQ_h\ve^{n+1}\right) + \left(B_he_{\theta}^{n+1} - B_he_{\theta}^n, B_he_{\theta}^{n+1}\right) \\\nonumber
			&\qquad+ \nu k\|\nab\ve^{n+1}\|^2_{\vL^2} + \mu k\|\nab e_{\theta}^{n+1}\|^2_{L^2}\\\nonumber
			&=\nu k\left(\nab \ve^{n+1},\nab\left(\vu(t_{n+1})-  \cQ_h\vu(t_{n+1})\right)\right) +\mu k\left(\nab e_{\theta}^{n+1},\nab\left(\theta(t_{n+1})-  B_h\theta(t_{n+1})\right)\right)  \\\nonumber
			&\qquad+ \nu \int_{t_n}^{t_{n+1}} \left(\nab(\vu(t_{n+1}) - \vu(s)),\nab \cQ_h\ve^{n+1}\right)\, ds \\\nonumber
			&\qquad+\mu \int_{t_n}^{t_{n+1}} \left(\nab(\theta(t_{n+1}) - \theta(s)),\nab B_h e_{\theta}^{n+1}\right)\, ds + k\left(\tilde{p}(t_{n+1}), \div \cQ_h\ve^{n+1}\right)\\\nonumber
			&\qquad + k \left(e_{\theta}^n{\bf e}_2, \cQ_h\ve^{n+1}\right) + \int_{t_{n}}^{t_{n+1}}\left((\theta(s) - \theta(t_n)){\bf e}_2, \cQ_h\ve^{n+1}\right)\, ds\\\nonumber
			&\qquad +\int_{t_n}^{t_{n+1}}\left[{b}(\vu(t_{n}),\vu(t_{n+1}),\cQ_h\ve^{n+1}) - {b}(\vu(s), \vu(s),\cQ_h\ve^{n+1}) \right] \, ds
			\\\nonumber
			&\qquad -k\left[{b}(\vu(t_{n}),\vu(t_{n+1}),\cQ_h\ve^{n+1}) - {b}(\vu_h^{n}, \vu_h^{n+1}, \cQ_h\ve^{n+1})  \right] \\\nonumber
			&\qquad + \int_{t_n}^{t_{n+1}}\left[\tilde{b}(\vu(t_{n+1}),\theta(t_{n+1}), B_h e^{n+1}_{\theta})- \tilde{b}(\vu(s),\theta(s),B_h e^{n+1}_{\theta})\right]\, ds    \\\nonumber
			&\qquad -k\left[\tilde{b}(\vu(t_{n+1}),\theta(t_{n+1}),B_h e^{n+1}_{\theta})\, ds - \tilde{b}(\vu_h^{n+1}, \theta_h^{n+1},B_h e^{n+1}_{\theta})  \right] \\\nonumber
			&\qquad+ \left(\int_{t_n}^{t_{n+1}}\left(\vG_1(\vu(s)) - \vG_1(\vu(t_n))\right)\, dW_1(s), \cQ_h\ve^{n+1}\right)\\\nonumber
			&\qquad+\left(\left(\vG_1(\vu(t_n))-\vG_1(\vu_h^n)\right)\Delta W_{1,n+1},\cQ_h\ve^{n+1}\right) \\\nonumber
			&\qquad+\left(\int_{t_n}^{t_{n+}} \left(G_2(\theta(s)) - G_2(\theta(t_n))\right)\, dW_2(s), B_h e^{n+1}_{\theta}\right)\\\nonumber
			&\qquad+\left(\left(G_2(\theta(t_n))- G_2(\theta_h^n)\right)\Delta W_{2,n+1},B_h e^{n+1}_{\theta}\right)\\\nonumber
			&:= I_1 + ... + I_{15}.
		\end{align}

		First, we notice that the left side of \eqref{eq_error} can be analyzed as follows:
		\begin{align*}
			\left(\cQ_h\ve^{n+1} - \cQ_h\ve^n, \cQ_h\ve^{n+1}\right)  &= \frac12\left[\|\cQ_h\ve^{n+1}\|^2_{\vL^2} - \|\cQ_h\ve^{n}\|^2_{\vL^2}\right] + \frac12 \|\cQ_h(\ve^{n+1} - \ve^n)\|^2_{\vL^2}.\\
			\left(B_h e_{\theta}^{n+1} - B_h e_{\theta}^n, B_h e_{\theta}^{n+1}\right)  &= \frac12\left[\|B_he_{\theta}^{n+1}\|^2_{L^2} - \|B_he_{\theta}^{n}\|^2_{L^2}\right] + \frac12 \|B_h(e_{\theta}^{n+1} - e_{\theta}^n)\|^2_{L^2}.
		\end{align*}
		Now, we estimate $I_1, ..., I_{15}$ as follows. First, using Cauchy-Schwarz's inequality and \eqref{Q_property_a}, and \eqref{B_property_b}, we have
		\begin{align*}
			I_1 + I_2 &\leq \frac{\nu k}{16}\|\nab \ve^{n+1}\|^2_{\vL^2} + Ck\|\nab(\vu(t_{n+1}) - \cQ_h\vu(t_{n+1}))\|^2_{\vL^2}\\\nonumber
			&\qquad+\frac{\mu k}{16}\|\nab e_\theta^{n+1}\|^2_{L^2} + Ck\|\nab(\theta(t_{n+1}) - B_h\theta(t_{n+1}))\|^2_{L^2}\\\nonumber
			&\leq \frac{\nu k}{16}\|\nab \ve^{n+1}\|^2_{\vL^2} + \frac{\mu k}{16}\|\nab e_{\theta}^{n+1}\|^2_{L^2} + Ckh^2\|\vA\vu(t_{n+1})\|^2_{\vL^2}  + Ckh^2\|\theta(t_{n+1})\|^2_{H^2}.
		\end{align*}
		
		Similarly, we estimate $I_3 + I_4$ as below.
		\begin{align*}
			I_3 + I_4 &\leq\frac{\nu k}{16}\|\nab \ve^{n+1}\|^2_{\vL^2} + \frac{\mu k}{16}\|\nab e_{\theta}^{n+1}\|^2_{L^2} \\\nonumber
			&\qquad+ C\int_{t_{n}}^{t_{n+1}} \left[\|\nab(\vu(t_{n+1}) - \vu(s))\|^2_{\vL^2} + \|\nab(\theta(t_{n+1}) - \theta(s))\|^2_{L^2}\right]\, ds.
		\end{align*}
		
		Next, using the fact that $\left(P_h \tilde{p}(t_{n+1}), \div\cQ_h\ve^{n+1}\right)=0$, we have 
		\begin{align*}
			I_5 &= k\left(\tilde{p}(t_{n+1}), \div \cQ_h\ve^{n+1}\right)\\\nonumber
			&=k\left(\tilde{p}(t_{n+1}) - P_h\tilde{p}(t_{n+1}), \div \cQ_h\ve^{n+1}\right)\\\nonumber
			&\leq \frac{\nu k}{16} \|\nab \ve^{n+1}\|^2_{\vL^2} + Ck \|\tilde{p}(t_{n+1}) - P_h\tilde{p}(t_{n+1})\|^2_{\vL^2}\\\nonumber
			&\leq \frac{\nu k}{16} \|\nab \ve^{n+1}\|^2_{\vL^2} + Ck h^2\|\nab \tilde{p}(t_{n+1})\|^2_{\vL^2},
		\end{align*} 
		where the last inequality was obtained by using \eqref{P_property_a}.
		
		Using the Cauchy-Schwarz inequality, we also get
		\begin{align*}
			I_6 + I_7 &\leq k \|e_{\theta}^n\|^2_{L^2} + \frac14 k \|\cQ_h \ve^{n+1}\|^2_{\vL^2} + \frac14 k \|\cQ_h\ve^{n+1}\|^2_{\vL^2} + \int_{t_{n}}^{t_{n+1}}\|\theta(t_{n+1}) - \theta(s)\|^2_{L^2}\, ds\\\nonumber
			&\leq 2k \|B_he_{\theta}^n\|^2_{L^2}  + Ckh^4\|\theta(t_n)\|^2_{H^2}+ \frac12 k \|\cQ_h\ve^{n+1}\|^2_{\vL^2} + \int_{t_{n}}^{t_{n+1}}\|\theta(t_{n+1}) - \theta(s)\|^2_{L^2}\, ds.
		\end{align*}
		
		Now, using the tri-linear property of ${b}(\cdot,\cdot,\cdot)$ and the Cauchy-Schwarz inequality, we also have
		\begin{align*}
			I_8 &= \int_{t_n}^{t_{n+1}}\left[{b}(\vu(t_{n}),\vu(t_{n+1}),\cQ_h\ve^{n+1}) - {b}(\vu(s), \vu(s),\cQ_h\ve^{n+1}) \right] \, ds\\\nonumber
			&=\int_{t_n}^{t_{n+1}}\left[{b}(\vu(t_{n}) - \vu(s),\vu(t_{n+1}),\cQ_h\ve^{n+1})+ {b}(\vu(s), \vu(t_{n+1})-\vu(s),\cQ_h\ve^{n+1}) \right] \, ds\\\nonumber
			&\leq \int_{t_{n}}^{t_{n+1}}\left[\|\vu(t_n) - \vu(s)\|_{\vL^4}\|\vu(t_{n+1})\|_{\vL^4}+ \|\vu(s)\|_{\vL^4}\|\vu(t_{n+1}) - \vu(s)\|_{\vL^4}\right]\|\nab\cQ_h\ve^{n+1}\|_{\vL^2}\, ds\\\nonumber
			&\leq \frac{\nu k}{16} \|\nab \ve^{n+1}\|^2_{\vL^2} + C\int_{t_n}^{t_{n+1}}\|\vu(t_{n}) - \vu(s)\|^2_{\vL^4}\|\vu(t_{n+1})\|^2_{\vL^4}\, ds \\\nonumber
			&\qquad+ C\int_{t_n}^{t_{n+1}}\|\vu(t_{n+1}) - \vu(s)\|^2_{\vL^4}\|\vu(s)\|^2_{\vL^4}\, ds.
		\end{align*}

		Next, we analyze and estimate $I_9$ as follows. First, we notice that 
		\begin{align*}
			\cQ_h\ve^{n+1} = \cQ_h\vu(t_{n+1}) - \vu_h^{n+1} = \ve^{n+1} - \left[\vu(t_{n+1}) - \cQ_h\vu(t_{n+1})\right].
		\end{align*}
		Therfore, with this and \eqref{b_property} we obtain
		\begin{align*}
			I_9 &= -k\left[{b}(\vu(t_{n}),\vu(t_{n+1}),\cQ_h\ve^{n+1}) - {b}(\vu_h^{n}, \vu_h^{n+1},\cQ_h\ve^{n+1})  \right]\\\nonumber
			&= -k{b}(\vu(t_{n}), \ve^{n+1},\cQ_h\ve^{n+1}) -k {b}(\ve^{n}, \vu_h^{n+1},\cQ_h\ve^{n+1})\\\nonumber
			&= -k {b}(\vu(t_{n}), \ve^{n+1},\ve^{n+1}) + k{b}(\vu(t_{n}), \ve^{n+1},\vu(t_{n+1}) - \cQ_h\vu(t_{n+1}))  \\\nonumber
			&\qquad+ k{b}(\ve^{n}, \ve^{n+1},\cQ_h\ve^{n+1}) - k{b}(\ve^{n}, \vu(t_{n+1}),\cQ_h\ve^{n+1})  \\\nonumber
			&= -k{b}(\vu(t_{n}), \ve^{n+1},\ve^{n+1}) + k{b}(\vu(t_{n}), \ve^{n+1},\vu(t_{n+1}) - \cQ_h\vu(t_{n+1}))  \\\nonumber
			&\qquad+ k{b}(\ve^{n}, \ve^{n+1},\ve^{n+1}) -k{b}(\ve^{n}, \ve^{n+1},\vu(t_{n+1})- \cQ_h \vu(t_{n+1})) \\\nonumber&\qquad\qquad- k{b}(\ve^{n}, \vu(t_{n+1}),\cQ_h\ve^{n+1})   \\\nonumber
			&= k{b}(\vu(t_{n}), \ve^{n+1},\vu(t_{n+1}) - \cQ_h\vu(t_{n+1}))  -k{b}(\ve^{n}, \ve^{n+1},\vu(t_{n+1})- \cQ_h \vu(t_{n+1})) \\\nonumber&\qquad\qquad- k{b}(\ve^{n}, \vu(t_{n+1}),\cQ_h\ve^{n+1})  \\\nonumber
			&:= I_{9,1} + I_{9,2} + I_{9,3}.
		\end{align*}
		
		Next, using the Ladyzhenskaya inequality and \eqref{Q_property_a}, we get
		\begin{align*}
			I_{9,1} &= k{b}(\vu(t_{n}), \ve^{n+1},\vu(t_{n+1}) - \cQ_h\vu(t_{n+1})) \\\nonumber
			&=k\left([\vu(t_{n})\cdot\nab]\ve^{n+1},\vu(t_{n+1}) - \cQ_h\vu(t_{n+1})\right)\\\nonumber
			&\leq k\|\vu(t_{n})\|_{\vL^4}\|\nab\ve^{n+1}\|_{\vL^2}\|\vu(t_{n+1}) - \cQ_h\vu(t_{n+1})\|_{\vL^4}\\\nonumber
			&\leq C k\|\vu(t_{n})\|_{\vH^1}\|\nab\ve^{n+1}\|_{\vL^2}\|\vu(t_{n+1}) - \cQ_h\vu(t_{n+1})\|^{\frac12}_{\vL^2}\|\nab(\vu(t_{n+1}) - \cQ_h\vu(t_{n+1}))\|^{\frac12}_{\vL^2}\\\nonumber
			&\leq \frac{\nu k}{16}\|\nab\ve^{n+1}\|^2_{\vL^2} + Ckh^3 \|\vu(t_{n})\|^2_{\vH^1}\|\vA\vu(t_{n+1})\|^2_{\vL^2}.
		\end{align*}
		
		Similarly, using the Ladyzhenskaya inequality and \eqref{Q_property_a}, we also obtain
		\begin{align*}
			I_{9,2} &= -k\left([\ve^{n}\cdot\nab]\ve^{n+1}, \vu(t_{n+1}) - \cQ_h\vu(t_{n+1})\right) - \frac{k}{2}\left([\div \ve^{n}]\ve^{n+1},\vu(t_{n+1}) - \cQ_h\vu(t_{n+1})\right)\\\nonumber
			&\leq k\|\ve^{n}\|_{\vL^4}\|\nab\ve^{n+1}\|_{\vL^2}\|\vu(t_{n+1}) - \cQ_h\vu(t_{n+1})\|_{\vL^4} \\\nonumber
			&\qquad+ \frac{k}{2}\|\nab\ve^{n}\|_{\vL^2}\|\ve^{n+1}\|_{\vL^4}\|\vu(t_{n+1}) - \cQ_h\vu(t_{n+1})\|_{\vL^4}\\\nonumber
			&\leq \frac{\nu k}{16}\left[\|\nab \ve^{n+1}\|^2_{\vL^2} + \|\nab \ve^{n}\|^2_{\vL^2}\right] + Ckh^3\|\ve^n\|^2_{\vL^4}\|\vA\vu(t_{n+1})\|^2_{\vL^2} \\\nonumber
			&\qquad+ Ckh^3\|\ve^{n+1}\|^2_{\vL^4}\|\vA\vu(t_{n+1})\|^2_{\vL^2}\\\nonumber
			&\leq \frac{\nu k}{16}\left[\|\nab \ve^{n+1}\|^2_{\vL^2} + \|\nab \ve^{n}\|^2_{\vL^2}\right] + Ckh^3\left[\|\vu(t_{n+1})\|^2_{\vH^1} + \|\vu(t_{n})\|^2_{\vH^1}\right]\|\vA\vu(t_{n+1})\|^2_{\vL^2} \\\nonumber
			&\qquad+ Ckh^2\left[\|\vu_h^{n+1}\|^2_{\vL^2} + \|\vu_h^{n}\|^2_{\vL^2}\right]\|\vA\vu(t_{n+1})\|^2_{\vL^2},
		\end{align*}
		where the last inequality was obtained by using the inverse inequality $\|\vu_h^{}\|_{\vL^4} \leq Ch^{-\frac12}\|\vu_h^{}\|_{\vL^2}$.
		
		Next, the term $I_{9,3}$ can be analyzed and bounded as follows.
		\begin{align*}
			I_{9,3} &= -k\left([\ve^{n}\cdot\nab]\vu(t_{n+1}),\cQ_h\ve^{n+1}\right) - \frac{k}{2}\left([\div \ve^{n}]\vu(t_{n+1}),\cQ_h\ve^{n+1}\right)\\\nonumber
			&:= I_{9,3a} + I_{9,3b}.
		\end{align*}
		
		Using the Ladyzhenskaya inequality and \eqref{Q_property_a}, the first term $I_{9,3a}$ can be estimated as below:
		\begin{align*}
			I_{9,3a} &= -k\left([\cQ_h\ve^{n}\cdot\nab]\vu(t_{n+1}),\cQ_h\ve^{n+1}\right) - k\left([(\vu(t_{n}) - \cQ_h\vu(t_{n}))\cdot\nab]\vu(t_{n+1}),\cQ_h\ve^{n+1}\right)\\\nonumber
			&\leq k\|\cQ_h\ve^{n}\|_{\vL^4}\|\cQ_h\ve^{n+1}\|_{\vL^4}\|\nab\vu(t_{n+1})\|_{\vL^2} + \frac{\nu k}{16}\|\nab \ve^{n+1}\|^2_{\vL^2} \\\nonumber
			&\qquad\qquad+ Ck\|\vu(t_{n}) - \cQ_h\vu(t_{n})\|^2_{\vL^4}\|\vu(t_{n+1})\|^2_{\vL^4}\\\nonumber
			&\leq C_e^2k\|\cQ_h\ve^{n}\|^{\frac12}_{\vL^2}\|\nab\cQ_h\ve^{n}\|^{\frac12}_{\vL^2}\|\cQ_h\ve^{n+1}\|^{\frac12}_{\vL^2}\|\nab\cQ_h\ve^{n+1}\|^{\frac12}_{\vL^2}\|\nab\vu(t_{n+1})\|_{\vL^2} + \frac{\nu k}{16}\|\nab \ve^{n+1}\|^2_{\vL^2} \\\nonumber
			&\qquad\qquad+ Ck\|\vu(t_{n}) - \cQ_h\vu(t_{n})\|^2_{\vL^4}\|\vu(t_{n+1})\|^2_{\vL^4}\\\nonumber
			&\leq \frac{C^2_ek}{\nu}\|\nab\vu(t_{n+1})\|^2_{\vL^2}\left[\|\cQ_h\ve^{n+1}\|^2_{\vL^2} + \|\cQ_h\ve^{n}\|^2_{\vL^2}\right] + \frac{\nu k}{16}\left[\|\nab \ve^{n+1}\|^2_{\vL^2} + \|\nab \ve^{n}\|^2_{\vL^2}\right] \\\nonumber
			&\qquad\qquad +\frac{\nu k}{16}\|\nab \ve^{n+1}\|^2_{\vL^2} + Ckh^3\|\vu(t_{n+1})\|^2_{\vH^1}\|\vA\vu(t_{n})\|^2_{\vL^2}\\\nonumber
			&\leq \frac{2C^2_ek}{\nu}\|\nab\vu(t_{n})\|^2_{\vL^2}\left[\|\cQ_h\ve^{n+1}\|^2_{\vL^2} + \|\cQ_h\ve^{n}\|^2_{\vL^2}\right] \\\nonumber
			&\qquad+\frac{2C^2_ek}{\nu}\|\nab(\vu(t_{n+1})-\vu(t_n))\|^2_{\vL^2}\left[\|\cQ_h\ve^{n+1}\|^2_{\vL^2} + \|\cQ_h\ve^{n}\|^2_{\vL^2}\right]  \\\nonumber
			&\qquad\qquad + \frac{\nu k}{16}\left[\|\nab \ve^{n+1}\|^2_{\vL^2} + \|\nab \ve^{n}\|^2_{\vL^2}\right] +\frac{\nu k}{16}\|\nab \ve^{n+1}\|^2_{\vL^2} + Ckh^3\|\vu(t_{n+1})\|^2_{\vH^1}\|\vA\vu(t_{n})\|^2_{\vL^2}.
		\end{align*}	
		
		Now, using H\"older's inequality and the Ladyzhenskaya inequality to estimate $I_{9,3b}$, we obtain
		\begin{align*}
			I_{9,3b} &= - \frac{k}{2}\left([\div \ve^{n}]\vu(t_{n+1}),\cQ_h\ve^{n+1}\right)\\\nonumber
			&\leq \frac{\nu k}{16} \|\nab \ve^{n}\|^2_{\vL^2} + \frac{4k}{\nu}\|\vu(t_{n+1})\cdot\cQ_h\ve^{n+1}\|^2_{\vL^2}\\\nonumber
			&\leq \frac{\nu k}{16} \|\nab \ve^{n}\|^2_{\vL^2} + \frac{4k}{\nu}\|\vu(t_{n+1})\|^2_{\vL^4}\|\cQ_h\ve^{n+1}\|^2_{\vL^4}\\\nonumber
			&\leq \frac{\nu k}{16} \|\nab \ve^{n}\|^2_{\vL^2} + \frac{4C_e^2k}{\nu}\|\vu(t_{n+1})\|^2_{\vH^1}\|\cQ_h\ve^{n+1}\|_{\vL^2}\|\nab\ve^{n+1}\|_{\vL^2}\\\nonumber
			&\leq \frac{\nu k}{16} \left[\|\nab \ve^{n}\|^2_{\vL^2} + \|\nab \ve^{n+1}\|^2_{\vL^2}\right] + \frac{64C_e^4k}{\nu^2}\|\vu(t_{n+1})\|^4_{\vH^1}\|\cQ_h\ve^{n+1}\|^2_{\vL^2}\\\nonumber&\leq \frac{\nu k}{16} \left[\|\nab \ve^{n}\|^2_{\vL^2} + \|\nab \ve^{n+1}\|^2_{\vL^2}\right] + \frac{512C_e^4k}{\nu^2}\|\vu(t_{n})\|^4_{\vH^1}\|\cQ_h\ve^{n+1}\|^2_{\vL^2} \\\nonumber
			&\qquad+ \frac{512C_e^4k}{\nu^2}\|\vu(t_{n+1}) - \vu(t_{n})\|^4_{\vH^1}\|\cQ_h\ve^{n+1}\|^2_{\vL^2}.
		\end{align*}
		
		Next, we estimate $I_{10}$ similarly as $I_8$. Using the Ladyzhenskaya inequality, we have
		\begin{align*}
			I_{10} &= \int_{t_n}^{t_{n+1}}\left[\tilde{b}(\vu(t_{n+1}),\theta(t_{n+1}),B_h e_{\theta}^{n+1}) - \tilde{b}(\vu(s), \theta(s), B_h e_{\theta}^{n+1}) \right] \, ds\\\nonumber
			&=\int_{t_n}^{t_{n+1}}\left[\tilde{b}(\vu(t_{n+1}) - \vu(s),\theta(t_{n+1}), B_h e_{\theta}^{n+1})+ \tilde{b}(\vu(s), \theta(t_{n+1})-\theta(s), B_h e_{\theta}^{n+1}) \right] \, ds\\\nonumber
			&\leq \int_{t_{n}}^{t_{n+1}}\left[\|\vu(t_{n+1}) - \vu(s)\|_{\vL^4}\|\theta(t_{n+1})\|_{L^4}+ \|\vu(s)\|_{\vL^4}\|\theta(t_{n+1}) - \theta(s)\|_{L^4}\right]\|\nab B_h e_{\theta}^{n+1}\|_{L^2}\, ds\\\nonumber
			&\leq \frac{\mu k}{16} \|\nab e_{\theta}^{n+1}\|^2_{L^2} + C\int_{t_n}^{t_{n+1}}\|\vu(t_{n+1}) - \vu(s)\|^2_{\vL^4}\|\theta(t_{n+1})\|^2_{L^4}\, ds \\\nonumber
			&\qquad+ C\int_{t_n}^{t_{n+1}}\|\theta(t_{n+1}) - \theta(s)\|^2_{L^4}\|\vu(s)\|^2_{\vL^4}\, ds\\\nonumber
			&\leq \frac{\mu k}{16} \|\nab e_{\theta}^{n+1}\|^2_{L^2} + C\int_{t_n}^{t_{n+1}}\|\vu(t_{n+1}) - \vu(s)\|^2_{\vH^1}\|\theta(t_{n+1})\|^2_{H^1}\, ds \\\nonumber
			&\qquad+ C\int_{t_n}^{t_{n+1}}\|\theta(t_{n+1}) - \theta(s)\|^2_{H^1}\|\vu(s)\|^2_{\vH^1}\, ds.
		\end{align*}
		
		Now, we are in the position to estimate $I_{11}$. Using the property \eqref{bb_property}, we obtain
		\begin{align*}
			I_{11} &= -k\left[\tilde{b}(\vu(t_{n+1}),\theta(t_{n+1}), B_h e_{\theta}^{n+1}) - \tilde{b}(\vu_h^{n+1}, \theta_h^{n+1}, B_h e_{\theta}^{n+1})  \right]\\\nonumber
			&= -k\tilde{b}(\vu(t_{n+1}), e_{\theta}^{n+1}, B_h e_{\theta}^{n+1}) -k \tilde{b}(\ve^{n+1}, \theta_h^{n+1}, B_h e_{\theta}^{n+1})\\\nonumber
			&= -k \tilde{b}(\vu(t_{n+1}), e_{\theta}^{n+1}, e_{\theta}^{n+1}) + k\tilde{b}(\vu(t_{n+1}), e_{\theta}^{n+1},\theta(t_{n+1}) - B_h\theta(t_{n+1}))  \\\nonumber
			&\qquad+ k\tilde{b}(\ve^{n+1}, e_{\theta}^{n+1},B_he_{\theta}^{n+1}) - k\tilde{b}(\ve^{n+1}, \theta(t_{n+1}), B_he_{\theta}^{n+1})  \\\nonumber
			&= -k\tilde{b}(\vu(t_{n+1}), e_{\theta}^{n+1}, e_{\theta}^{n+1}) + k\tilde{b}(\vu(t_{n+1}), e_{\theta}^{n+1},\theta(t_{n+1}) - B_h\theta(t_{n+1}))  \\\nonumber
			&\qquad+ k\tilde{b}(\ve^{n+1}, e_{\theta}^{n+1}, e_{\theta}^{n+1}) -k\tilde{b}(\ve^{n+1}, e_{\theta}^{n+1},\theta(t_{n+1})- B_h \theta(t_{n+1})) \\\nonumber&\qquad\qquad- k\tilde{b}(\ve^{n+1}, \theta(t_{n+1}), B_h e_{\theta}^{n+1})   \\\nonumber
			&= k\tilde{b}(\vu(t_{n+1}), e_{\theta}^{n+1},\theta(t_{n+1}) - B_h\theta(t_{n+1}))   -k\tilde{b}(\ve^{n+1}, e_{\theta}^{n+1},\theta(t_{n+1})- B_h \theta(t_{n+1})) \\\nonumber&\qquad\qquad- k\tilde{b}(\ve^{n+1}, \theta(t_{n+1}), B_h e_{\theta}^{n+1})   \\\nonumber
			&:= I_{11,1} + I_{11,2} + I_{11,3}.
		\end{align*}
		
		Next, using the Ladyzhenskaya inequality and \eqref{B_property_b}, we obtain
		\begin{align*}
			I_{11,1} &= k\tilde{b}(\vu(t_{n+1}), e_{\theta}^{n+1},\theta(t_{n+1}) - B_h\theta(t_{n+1}))   \\\nonumber
			&=k\left(\vu(t_{n+1})\cdot\nab e_{\theta}^{n+1},\theta(t_{n+1}) - B_h\theta(t_{n+1})\right)\\\nonumber
			&\leq k\|\vu(t_{n+1})\|_{\vL^4}\|\nab e_{\theta}^{n+1}\|_{L^2}\|\theta(t_{n+1}) - B_h\theta(t_{n+1})\|_{L^4}\\\nonumber
			&\leq C k\|\vu(t_{n+1})\|_{\vH^1}\|\nab e_{\theta}^{n+1}\|_{L^2}\|\theta(t_{n+1}) - B_h\theta(t_{n+1})\|^{\frac12}_{L^2}\|\nab(\theta(t_{n+1}) - B_h\theta(t_{n+1}))\|^{\frac12}_{L^2}\\\nonumber
			&\leq \frac{\mu k}{16}\|\nab e_{\theta}^{n+1}\|^2_{L^2} + Ckh^3 \|\vu(t_{n+1})\|^2_{\vH^1}\|\theta(t_{n+1})\|^2_{H^2}.
		\end{align*}
		
		To estimate $I_{11,2}$, we also using the Ladyzhenskaya inequality, \eqref{B_property_b} and the inverse inequality $\|\vu_h^{}\|_{\vL^4} \leq Ch^{-\frac12}\|\vu_h^{}\|_{\vL^2}$ as follows:
		\begin{align*}
			I_{11,2} &= -k\left(\ve^{n+1}\cdot\nab e_{\theta}^{n+1}, \theta(t_{n+1}) - B_h\theta(t_{n+1})\right) - \frac{k}{2}\left([\div \ve^{n+1}] e_{\theta}^{n+1},\theta(t_{n+1}) - B_h\theta(t_{n+1})\right)\\\nonumber
			&\leq k\|\ve^{n+1}\|_{\vL^4}\|\nab e_{\theta}^{n+1}\|_{L^2}\|\theta(t_{n+1}) - B_h\theta(t_{n+1})\|_{L^4} \\\nonumber
			&\qquad+ \frac{k}{2}\|\nab\ve^{n+1}\|_{\vL^2}\|e_{\theta}^{n+1}\|_{L^4}\|\theta(t_{n+1}) - B_h\theta(t_{n+1})\|_{L^4}\\\nonumber
			&\leq \frac{\mu k}{16}\|\nab e_{\theta}^{n+1}\|^2_{L^2} + \frac{\nu k}{16}\|\nab \ve^{n+1}\|^2_{\vL^2}+ Ckh^3\|\ve^{n+1}\|^2_{\vL^4}\|\theta(t_{n+1})\|^2_{H^2} + Ckh^3\|e_{\theta}^{n+1}\|^2_{L^4}\|\theta(t_{n+1})\|^2_{H^2}\\\nonumber
			&\leq \frac{\mu k}{16}\|\nab e_{\theta}^{n+1}\|^2_{L^2} + \frac{\nu k}{16}\|\nab \ve^{n+1}\|^2_{\vL^2} + Ckh^3\left[\|\theta(t_{n+1})\|^2_{H^1} + \|\vu(t_{n+1})\|^2_{\vH^1}\right]\|\theta(t_{n+1})\|^2_{H^2} \\\nonumber
			&\qquad+ Ckh^2\left[\|\theta_h^{n+1}\|^2_{L^2} + \|\vu_h^{n+1}\|^2_{\vL^2}\right]\|\theta(t_{n+1})\|^2_{H^2}.
		\end{align*}
		
		Next, we are in the position to estimate $I_{11,3}$, which can be analyzed as follows:
		\begin{align*}
			I_{11,3} &= -k\left(\ve^{n+1}\cdot\nab\theta(t_{n+1}), B_h e_{\theta}^{n+1}\right) - \frac{k}{2}\left([\div \ve^{n+1}]\theta(t_{n+1}), B_h e_{\theta}^{n+1}\right)\\\nonumber
			&:= I_{11,3a} + I_{11,3b}.
		\end{align*}
		
		Using the Ladyzhenskaya inequality and \eqref{B_property_b}, we have
		\begin{align*}
			I_{11,3a} &= -k\left(\cQ_h\ve^{n+1}\cdot\nab\theta(t_{n+1}), B_he_{\theta}^{n+1}\right) - k\left((\vu(t_{n+1}) - \cQ_h\vu(t_{n+1}))\cdot\nab\theta(t_{n+1}), B_he_{\theta}^{n+1}\right)\\\nonumber
			&\leq k\|\cQ_h\ve^{n+1}\|_{\vL^4}\|B_he_{\theta}^{n+1}\|_{L^4}\|\nab\theta(t_{n+1})\|_{L^2} + \frac{\mu k}{16}\|\nab e_{\theta}^{n+1}\|^2_{L^2} \\\nonumber
			&\qquad\qquad+ Ck\|\vu(t_{n+1}) - \cQ_h\vu(t_{n+1})\|^2_{\vL^4}\|\theta(t_{n+1})\|^2_{L^4}\\\nonumber
			&\leq C_e^2k\|\cQ_h\ve^{n+1}\|^{\frac12}_{\vL^2}\|\nab\cQ_h\ve^{n+1}\|^{\frac12}_{\vL^2}\|B_he_{\theta}^{n+1}\|^{\frac12}_{L^2}\|\nab B_he_{\theta}^{n+1}\|^{\frac12}_{L^2}\|\nab\theta(t_{n+1})\|_{L^2} \\\nonumber
			&\qquad\qquad+ \frac{\mu k}{16}\|\nab e_{\theta}^{n+1}\|^2_{L^2} + Ck\|\vu(t_{n+1}) - \cQ_h\vu(t_{n+1})\|^2_{\vL^4}\|\theta(t_{n+1})\|^2_{L^4}\\\nonumber
			&\leq \frac{C^2_ek}{\nu}\|\nab\theta(t_{n+1})\|^2_{L^2}\left[\|B_he_{\theta}^{n+1}\|^2_{L^2} + \|\cQ_h\ve^{n}\|^2_{\vL^2}\right] + \frac{\mu k}{16}\|\nab e_{\theta}^{n+1}\|^2_{L^2} + \frac{\nu k}{16}\|\nab \ve^{n+1}\|^2_{\vL^2}\\\nonumber
			&\qquad\qquad +\frac{\mu k}{16}\|\nab e_{\theta}^{n+1}\|^2_{L^2} + Ckh^3\|\theta(t_{n+1})\|^2_{H^1}\|\vA\vu(t_{n+1})\|^2_{\vL^2}\\\nonumber
			&\leq \frac{2C^2_ek}{\nu}\|\nab\theta(t_{n})\|^2_{L^2}\left[\|B_he_{\theta}^{n+1}\|^2_{L^2} + \|\cQ_h\ve^{n}\|^2_{\vL^2}\right] \\\nonumber
			&\qquad+\frac{2C^2_ek}{\nu}\|\nab(\theta(t_{n+1})- \theta(t_{n}))\|^2_{L^2}\left[\|B_he_{\theta}^{n+1}\|^2_{L^2} + \|\cQ_h\ve^{n}\|^2_{\vL^2}\right]  + \frac{\nu k}{16}\|\nab \ve^{n+1}\|^2_{\vL^2}\\\nonumber
			&\qquad\qquad +\frac{\mu k}{8}\|\nab e_{\theta}^{n+1}\|^2_{L^2} + Ckh^3\|\theta(t_{n+1})\|^2_{H^1}\|\vA\vu(t_{n+1})\|^2_{\vL^2}.
		\end{align*}	
		
		Now, using the Ladyzhenskaya inequality to estimate $I_{11,3b}$, we obtain
		\begin{align*}
			I_{11,3b} &= - \frac{k}{2}\left([\div \ve^{n+1}]\theta(t_{n+1}), B_h e_{\theta}^{n+1}\right)\\\nonumber
			&\leq \frac{\nu k}{16} \|\nab \ve^{n+1}\|^2_{\vL^2} + \frac{4k}{\nu}\|\theta(t_{n+1})\cdot B_he_{\theta}^{n+1}\|^2_{L^2}\\\nonumber
			&\leq \frac{\nu k}{16} \|\nab \ve^{n+1}\|^2_{\vL^2} + \frac{4k}{\nu}\|\theta(t_{n+1})\|^2_{L^4}\| B_h e_{\theta}^{n+1}\|^2_{L^4}\\\nonumber
			&\leq \frac{\nu k}{16} \|\nab \ve^{n+1}\|^2_{\vL^2} + \frac{4C_e^2k}{\nu}\|\theta(t_{n+1})\|^2_{H^1}\|B_he_{\theta}^{n+1}\|_{L^2}\|\nab e_{\theta}^{n+1}\|_{L^2}\\\nonumber
			&\leq \frac{\nu k}{16} \|\nab \ve^{n+1}\|^2_{\vL^2} + \frac{\mu k}{16}\|\nab e_{\theta}^{n+1}\|^2_{L^2} + \frac{64C_e^4k}{\nu^2}\|\theta(t_{n+1})\|^4_{H^1}\|B_he_{\theta}^{n+1}\|^2_{L^2}\\\nonumber
			&\leq \frac{\nu k}{16} \|\nab \ve^{n+1}\|^2_{\vL^2} + \frac{\mu k}{16}\|\nab e_{\theta}^{n+1}\|^2_{L^2} + \frac{512C_e^4k}{\nu^2}\|\theta(t_{n})\|^4_{H^1}\|B_he_{\theta}^{n+1}\|^2_{L^2} \\\nonumber
			&\qquad+ \frac{512C_e^4k}{\nu^2}\|\theta(t_{n+1}) - \theta(t_{n})\|^4_{H^1}\|B_he_{\theta}^{n+1}\|^2_{L^2}.
		\end{align*}
		
		Now, we are dealing with the noise terms. First, we have
		\begin{align*}
			I_{12} + I_{14} + I_{13} + I_{15}&=  \left(\int_{t_n}^{t_{n+1}}\left(\vG_1(\vu(s)) - \vG_1(\vu(t_n))\right)\, dW_1(s), \cQ_h\ve^{n+1} - \cQ_h\ve^n\right)\\\nonumber
			&\qquad+\left(\int_{t_n}^{t_{n+}} \left(G_2(\theta(s)) - G_2(\theta(t_n))\right)\, dW_2(s), B_h e^{n+1}_{\theta} - B_h e^{n}_{\theta}\right)\\\nonumber
			&\qquad+  \left(\int_{t_n}^{t_{n+1}}\left(\vG_1(\vu(s)) - \vG_1(\vu(t_n))\right)\, dW_1(s),  \cQ_h\ve^n\right)\\\nonumber
			&\qquad+\left(\int_{t_n}^{t_{n+}} \left(G_2(\theta(s)) - G_2(\theta(t_n))\right)\, dW_2(s),  B_h e^{n}_{\theta}\right) \\\nonumber
			&\qquad+\left(\left(\vG_1(\vu(t_n))-\vG_1(\vu_h^n)\right)\Delta W_{1,n+1},\cQ_h\ve^{n+1} - \cQ_h\ve^n\right) \\\nonumber
			&\qquad+\left(\left(G_2(\theta(t_n))- G_2(\theta_h^n)\right)\Delta W_{2,n+1}, B_he_{\theta}^{n+1} - B_h e_{\theta}^{n}\right)  \\\nonumber
			&\qquad+\left(\left(\vG_1(\vu(t_n))-\vG_1(\vu_h^n)\right)\Delta W_{1,n+1},\cQ_h\ve^{n}\right)  \\\nonumber
			&\qquad+\left(\left(G_2(\theta(t_n))-G_2(\theta_h^n)\right)\Delta W_{2,n+1}, B_h e_{\theta}^{n+1}\right), 
		\end{align*}
		which, together with the Young inequality, and the assumptions (B3) and (B4), and \eqref{Q_property_a}, and \eqref{B_property_b} implies that
		\begin{align*}
			I_{12} + I_{14} + I_{13} + I_{15}&\leq \frac14\left[\|\cQ_h\ve^{n+1} - \cQ_h\ve^n\|^2_{\vL^2} + \|B_h e_{\theta}^{n+1} - B_he_{\theta}^n\|^2_{L^2}\right] \\\nonumber&\qquad +2\left\|\int_{t_n}^{t_{n+1}}\left(\vG_1(\vu(s)) - \vG_1(\vu(t_n))\right)\, dW_1(s)\right\|^2_{\vL^2}\\\nonumber
			&\qquad+2\left\|\int_{t_n}^{t_{n+}} \left(G_2(\theta(s)) - G_2(\theta(t_n))\right)\, dW_2(s)\right\|^2_{L^2}\\\nonumber
			&\qquad + 4C_{\vG_1}\|\cQ_h\ve^n\|^2_{\vL^2}|\Delta W_{1,n+1}|^2 + 4C_{G_2}\|B_he_{\theta}^n\|^2_{L^2}|\Delta W_{2,n+1}|^2\\\nonumber
			&\qquad+ Ch^4\left[\|\vA\vu(t_n)\|^2_{\vL^2}|\Delta W_{1,n+1}|^2 + \|\theta(t_n)\|^2_{H^2}|\Delta W_{2,n+1}|^2\right]\\\nonumber
			&\qquad+  \left(\int_{t_n}^{t_{n+1}}\left(\vG_1(\vu(s)) - \vG_1(\vu(t_n))\right)\, dW_{1}(s),  \cQ_h\ve^n\right)\\\nonumber
			&\qquad+\left(\int_{t_n}^{t_{n+}} \left(G_2(\theta(s)) - G_2(\theta(t_n))\right)\, dW_2(s),  B_h e^{n}_{\theta}\right)  \\\nonumber
			&\qquad+\left(\left(\vG_1(\vu(t_n))-\vG_1(\vu_h^n)\right)\Delta W_{1,n+1},\cQ_h\ve^{n}\right)  \\\nonumber
			&\qquad+\left(\left(G_2(\theta(t_n))-G_2(\theta_h^n)\right)\Delta W_{2,n+1}, B_h e_{\theta}^{n}\right).
		\end{align*}
		Now, substituting all the estimates from $I_1, ..., I_{15}$ into the right-hand side of \eqref{eq_error}, and absorbing the like terms from the left-hand side to the right-hand side, then multiplying the result by the indicator function $\mathbf{1}_{\Omega_{\rho, n}}$ and also using the fact that $\mathbf{1}_{\Omega_{\rho, n}} \geq \mathbf{1}_{\Omega_{\rho, n+1}}$ we arrive at
		\begin{align}\label{eq3.22}
			&	\frac12\mathbf{1}_{\Ome,n+1}\|\cQ_h\ve^{n+1}\|^2_{\vL^2} - \frac12\mathbf{1}_{\Omega_{\rho, n}}\|\cQ_h\ve^n\|^2_{\vL^2} +  	\frac12\mathbf{1}_{\Ome,n+1}\|B_he_{\theta}^{n+1}\|^2_{L^2} - \frac12\mathbf{1}_{\Omega_{\rho, n}}\|B_he_{\theta}^n\|^2_{L^2}\\\nonumber
			&\qquad+ \frac{\nu k}{16}\mathbf{1}_{\Omega_{\rho, n}}\|\nab\ve^{n+1}\|^2_{\vL^2} + \frac{\mu k}{2}\mathbf{1}_{\Omega_{\rho, n}}\|\nab e_{\theta}^{n+1}\|^2_{L^2} \\\nonumber
			&\qquad+ \frac14\mathbf{1}_{\Omega_{\rho, n}}\left[\|\cQ_h\ve^{n+1} - \cQ_h\ve^n\|^2_{\vL^2} + \|B_h e_{\theta}^{n+1} - B_he_{\theta}^n\|^2_{L^2}\right]\\\nonumber
			&\leq \Biggl\{Ckh^2\|\vA\vu(t_{n+1})\|^2_{\vL^2}  + Ckh^2\|\theta(t_{n+1})\|^2_{H^2} + Ck h^2\|\nab \tilde{p}(t_{n+1})\|^2_{\vL^2}
			\\\nonumber
			&\qquad+  Ckh^3 \|\vu(t_{n})\|^2_{\vH^1}\|\vA\vu(t_{n+1})\|^2_{\vL^2} \\\nonumber
			&\qquad+ Ckh^3\left[\|\vu(t_{n+1})\|^2_{\vH^1} + \|\vu(t_{n})\|^2_{\vH^1}\right]\|\vA\vu(t_{n+1})\|^2_{\vL^2} \\\nonumber
			&\qquad+ Ckh^2\left[\|\vu_h^{n+1}\|^2_{\vL^2} + \|\vu_h^{n}\|^2_{\vL^2}\right]\|\vA\vu(t_{n+1})\|^2_{\vL^2}
			+ Ckh^3 \|\vu(t_{n+1})\|^2_{\vH^1}\|\theta(t_{n+1})\|^2_{H^2} \\\nonumber&\qquad+ Ckh^3\left[\|\theta(t_{n+1})\|^2_{H^1} + \|\vu(t_{n+1})\|^2_{\vH^1}\right]\|\theta(t_{n+1})\|^2_{H^2} \\\nonumber&\qquad+ Ckh^2\left[\|\theta_h^{n+1}\|^2_{L^2} + \|\vu_h^{n+1}\|^2_{\vL^2}\right]\|\theta(t_{n+1})\|^2_{H^2}
			\\\nonumber
			&\qquad+ C\int_{t_{n}}^{t_{n+1}} \left[\|\nab(\vu(t_{n+1}) - \vu(s))\|^2_{\vL^2} + \|\nab(\theta(t_{n+1}) - \theta(s))\|^2_{L^2}\right]\, ds  \\\nonumber
			&\qquad + \int_{t_{n}}^{t_{n+1}}\|\theta(t_{n+1}) - \theta(s)\|^2_{L^2}\, ds + \frac{512C_e^4k}{\nu^2}\|\vu(t_{n+1}) - \vu(t_{n})\|^4_{\vH^1}\|\cQ_h\ve^{n+1}\|^2_{\vL^2} \\\nonumber
			&\qquad+  C\int_{t_n}^{t_{n+1}}\|\vu(t_{n}) - \vu(s)\|^2_{\vL^4}\|\vu(t_{n+1})\|^2_{\vL^4}\, ds \\\nonumber
			&\qquad+ C\int_{t_n}^{t_{n+1}}\|\vu(t_{n+1}) - \vu(s)\|^2_{\vL^4}\|\vu(s)\|^2_{\vL^4}\, ds
			\\\nonumber
			&\qquad+\frac{2C^2_ek}{\nu}\|\nab(\vu(t_{n+1})-\vu(t_n))\|^2_{\vL^2}\left[\|\cQ_h\ve^{n+1}\|^2_{\vL^2} + \|\cQ_h\ve^{n}\|^2_{\vL^2}\right] 
			\\\nonumber
			&\qquad+ C\int_{t_n}^{t_{n+1}}\|\vu(t_{n+1}) - \vu(s)\|^2_{\vH^1}\|\theta(t_{n+1})\|^2_{H^1}\, ds \\\nonumber
			&\qquad+ C\int_{t_n}^{t_{n+1}}\|\theta(t_{n+1}) - \theta(s)\|^2_{H^1}\|\vu(s)\|^2_{\vH^1}\, ds	
			\\\nonumber
			&\qquad+\frac{2C^2_ek}{\nu}\|\nab(\theta(t_{n+1})- \theta(t_{n}))\|^2_{L^2}\left[\|B_he_{\theta}^{n+1}\|^2_{L^2} + \|\cQ_h\ve^{n}\|^2_{\vL^2}\right]  \\\nonumber
			&\qquad+ \frac{512C_e^4k}{\nu^2}\|\theta(t_{n+1}) - \theta(t_{n})\|^4_{H^1}\|B_he_{\theta}^{n+1}\|^2_{L^2}\\\nonumber&\qquad +2\left\|\int_{t_n}^{t_{n+1}}\left(\vG_1(\vu(s)) - \vG_1(\vu(t_n))\right)\, dW_1(s)\right\|^2_{\vL^2}\\\nonumber
			&\qquad+2\left\|\int_{t_n}^{t_{n+}} \left(G_2(\theta(s)) - G_2(\theta(t_n))\right)\, dW_@(s)\right\|^2_{L^2} 
			\\\nonumber
			&\qquad+  \mathbf{1}_{\Omega_{\rho, n}}\left(\int_{t_n}^{t_{n+1}}\left(\vG_1(\vu(s)) - \vG_1(\vu(t_n))\right)\, dW_1(s),  \cQ_h\ve^n\right)\\\nonumber
			&\qquad+  \mathbf{1}_{\Omega_{\rho, n}}\left(\int_{t_n}^{t_{n+}} \left(G_2(\theta(s)) - G_2(\theta(t_n))\right)\, dW_2(s),  B_h e^{n}_{\theta}\right)  \\\nonumber
			&\qquad+  \mathbf{1}_{\Omega_{\rho, n}}\left(\left(\vG_1(\vu(t_n))-\vG_1(\vu_h^n)\right)\Delta W_{1,n+1},\cQ_h\ve^{n}\right)  \\\nonumber
			&\qquad+  \mathbf{1}_{\Omega_{\rho, n}}\left(\left(G_2(\theta(t_n))-G_2(\theta_h^n)\right)\Delta W_{2,n+1}, B_h e_{\theta}^{n}\right) \\\nonumber
			&\qquad+ Ch^4 \left[\|\vA \vu(t_n)\|^2_{\vL^2}|\Delta W_{1,n+1}|^2 + \|\theta(t_n)\|^2_{H^2}|\Delta W_{2,n+1}|^2\right] \Biggr\}
			\\\nonumber&\qquad+ \Biggl\{k \mathbf{1}_{\Omega_{\rho, n}} \|B_he_{\theta}^n\|^2_{L^2} + \frac12 k  \mathbf{1}_{\Omega_{\rho, n}}\|\cQ_h \ve^{n+1}\|^2_{\vL^2} \\\nonumber
			&\qquad+\frac{2C^2_ek}{\nu} \mathbf{1}_{\Omega_{\rho, n}}\|\nab\vu(t_{n})\|^2_{\vL^2}\left[\|\cQ_h\ve^{n+1}\|^2_{\vL^2} + \|\cQ_h\ve^{n}\|^2_{\vL^2}\right] \\\nonumber
			&\qquad+ \frac{512C_e^4k}{\nu^2} \mathbf{1}_{\Omega_{\rho, n}}\|\vu(t_{n})\|^4_{\vH^1}\|\cQ_h\ve^{n+1}\|^2_{\vL^2} + \frac{512C_e^4k}{\nu^2} \mathbf{1}_{\Omega_{\rho, n}}\|\theta(t_{n})\|^4_{H^1}\|B_he_{\theta}^{n+1}\|^2_{L^2} \\\nonumber
			&\qquad +\frac{2C^2_ek}{\nu} \mathbf{1}_{\Omega_{\rho, n}}\|\nab\theta(t_{n})\|^2_{L^2}\left[\|B_he_{\theta}^{n+1}\|^2_{L^2} + \|\cQ_h\ve^{n}\|^2_{\vL^2}\right]\\\nonumber
			&\qquad + 4C_{\vG_1} \mathbf{1}_{\Omega_{\rho, n}}\|\cQ_h\ve^n\|^2_{\vL^2}|\Delta W_{1,n+1}|^2 + 4C_{G_2} \mathbf{1}_{\Omega_{\rho, n}}\|B_h e_{\theta}^{n}\|^2_{L^2}|\Delta W_{2,n+1}|^2\Biggr\}\\\nonumber
			&:= F_{n+1} + \Biggl\{k \mathbf{1}_{\Omega_{\rho, n}} \|B_he_{\theta}^n\|^2_{L^2} + \frac12 k  \mathbf{1}_{\Omega_{\rho, n}}\|\cQ_h \ve^{n+1}\|^2_{\vL^2} \\\nonumber
			&\qquad+\frac{2C^2_ek}{\nu} \mathbf{1}_{\Omega_{\rho, n}}\|\nab\vu(t_{n})\|^2_{\vL^2}\left[\|\cQ_h\ve^{n+1}\|^2_{\vL^2} + \|\cQ_h\ve^{n}\|^2_{\vL^2}\right] \\\nonumber
			&\qquad+ \frac{512C_e^4k}{\nu^2} \mathbf{1}_{\Omega_{\rho, n}}\|\vu(t_{n})\|^4_{\vH^1}\|\cQ_h\ve^{n+1}\|^2_{\vL^2} + \frac{512C_e^4k}{\nu^2} \mathbf{1}_{\Omega_{\rho, n}}\|\theta(t_{n})\|^4_{H^1}\|B_he_{\theta}^{n+1}\|^2_{L^2} \\\nonumber
			&\qquad +\frac{2C^2_ek}{\nu} \mathbf{1}_{\Omega_{\rho, n}}\|\nab\theta(t_{n})\|^2_{L^2}\left[\|B_he_{\theta}^{n+1}\|^2_{L^2} + \|\cQ_h\ve^{n}\|^2_{\vL^2}\right]\\\nonumber
			&\qquad + 4C_{\vG_1} \mathbf{1}_{\Omega_{\rho, n}}\|\cQ_h\ve^n\|^2_{\vL^2}|\Delta W_{1,n+1}|^2 + 4C_{G_2} \mathbf{1}_{\Omega_{\rho, n}}\|B_h e_{\theta}^{n}\|^2_{L^2}|\Delta W_{2,n+1}|^2\Biggr\}.
		\end{align}
		
		Next, taking the expectation and applying the summation $\sum_{n = 0}^{\ell}$ for any $0\leq \ell \leq M-1$ to \eqref{eq3.22}, we then obtain
		\begin{align}\label{eq3.23}
			&	\frac12\mE\left[\mathbf{1}_{\Ome,\ell+1}\|\cQ_h\ve^{\ell+1}\|^2_{\vL^2} +  	\frac12\mathbf{1}_{\Ome,\ell+1}\|B_he_{\theta}^{\ell+1}\|^2_{L^2}\right] \\\nonumber
			&\qquad+ \frac{\nu k}{16}\sum_{n =0}^{\ell}\mE\left[\mathbf{1}_{\Omega_{\rho, n}}\|\nab\ve^{n+1}\|^2_{\vL^2} + \frac{\mu k}{2}\mathbf{1}_{\Omega_{\rho, n}}\|\nab e_{\theta}^{n+1}\|^2_{L^2}\right] \\\nonumber
			&\leq \sum_{n = 0}^{\ell}\mE\left[F_{n+1}\right] + \sum_{n = 0}^{\ell}k \left(\frac32 + \tilde{C}\rho + 4C_{\vG_1}+ 4 C_{G_2}\right)\mE\left[\mathbf{1}_{\Omega_{\rho, n}}\left(\|\cQ_h\ve^{n+1}\|^2_{\vL^2} + \|\cQ_h\ve^{n}\|^2_{\vL^2} \right.\right.\\\nonumber
			&\qquad\left.\left.+ \|B_he_{\theta}^{n+1}\|^2_{L^2} + \|B_he_{\theta}^{n}\|^2_{L^2}\right)\right],
		\end{align}
		where $\tilde{C} = \frac{4C_e^2}{\nu} + \frac{1024 C_e^4}{\nu^2}$.
		
		We note that \eqref{eq3.23} gives us the applicable form of the discrete Gronwall inequality. So, it is left to estimate the first term, $\sum_{n = 0}^{\ell}\mE\left[F_{n+1}\right]$, on the right-hand side of \eqref{eq3.23}. First of all, we can group it into three groups $Z_1, Z_2$, and $Z_3$ as shown below. Each group can be estimated by using similar techniques. 
		
		To the end, 
		first we estimate $Z_1$.  Using Lemma \ref{stability_pdes}, Lemma \ref{stability_FEMs} and Lemma \ref{lemma2.2}, we obtain
		\begin{align*}
			Z_1 &= \sum_{n=0}^{\ell}\mE\biggl[Ckh^2\|\vA\vu(t_{n+1})\|^2_{\vL^2}  + Ckh^2\|\theta(t_{n+1})\|^2_{H^2} + Ck h^2\|\nab \tilde{p}(t_{n+1})\|^2_{\vL^2}
			\\\nonumber
			&\qquad+  Ckh^3 \|\vu(t_{n})\|^2_{\vH^1}\|\vA\vu(t_{n+1})\|^2_{\vL^2} + Ckh^3 \|\vu(t_{n+1})\|^2_{\vH^1}\|\theta(t_{n+1})\|^2_{H^2} \\\nonumber
			&\qquad+ Ckh^3\left[\|\vu(t_{n+1})\|^2_{\vH^1} + \|\vu(t_{n})\|^2_{\vH^1}\right]\|\vA\vu(t_{n+1})\|^2_{\vL^2} \\\nonumber
			&\qquad+ Ckh^2\left[\|\vu_h^{n+1}\|^2_{\vL^2} + \|\vu_h^{n}\|^2_{\vL^2}\right]\|\vA\vu(t_{n+1})\|^2_{\vL^2}
			\\\nonumber&\qquad+ Ckh^3\left[\|\theta(t_{n+1})\|^2_{H^1} + \|\vu(t_{n+1})\|^2_{\vH^1}\right]\|\theta(t_{n+1})\|^2_{H^2} \\\nonumber&\qquad+ Ckh^2\left[\|\theta_h^{n+1}\|^2_{L^2} + \|\vu_h^{n+1}\|^2_{\vL^2}\right]\|\theta(t_{n+1})\|^2_{H^2} \biggr]\\\nonumber
			&\leq Ch^2 + Ch^2k\sum_{n = 0}^{\ell}\mE\left[\|\nab \tilde{p}(t_{n+1})\|^2_{L^2}\right] \\\nonumber
			&\qquad+ Ch^3 \left(\mE\left[\sup_{s\in [0,T]}\|\vu(s)\|^4_{\vH^1}\right]\right)^{\frac12} \left(\mE\left[\left(k\sum_{n=0}^{\ell}\|\vA\vu(t_{n+1})\|^2_{\vL^2}\right)^2\right]\right)^{\frac12}\\\nonumber
			&\qquad+ Ch^3 \left(\mE\left[\sup_{s\in [0,T]}\|\vu(s)\|^4_{\vH^1}\right]\right)^{\frac12} \left(\mE\left[\left(k\sum_{n=0}^{\ell}\|\theta(t_{n+1})\|^2_{H^2}\right)^2\right]\right)^{\frac12}\\\nonumber
			&\qquad+ Ch^2 \left(\mE\left[\max_{1\leq n \leq M}\|\vu^n_h\|^4_{\vL^2}\right]\right)^{\frac12} \left(\mE\left[\left(k\sum_{n=0}^{\ell}\|\vA\vu(t_{n+1})\|^2_{\vL^2}\right)^2\right]\right)^{\frac12}\\\nonumber
			&\qquad+ Ch^3 \left(\mE\left[\sup_{s\in [0,T]}(\|\theta(s)\|^4_{H^1}+\|\vu(s)\|^4_{\vH^1})\right]\right)^{\frac12} \left(\mE\left[\left(k\sum_{n=0}^{\ell}\|\theta(t_{n+1})\|^2_{H^2}\right)^2\right]\right)^{\frac12}\\\nonumber
			&\qquad+ Ch^2 \left(\mE\left[\max_{1\leq n\leq M}(\|\theta^n_h\|^4_{L^2}+\|\vu^n_h\|^4_{\vL^2})\right]\right)^{\frac12} \left(\mE\left[\left(k\sum_{n=0}^{\ell}\|\theta(t_{n+1})\|^2_{H^2}\right)^2\right]\right)^{\frac12}\\\nonumber
			&\leq Ch^2,
		\end{align*}
		where the last inequality is obtained by using Lemma \ref{lemma2.2} (d) and the fact that
		\begin{align*}
			Ch^2k\sum_{n = 0}^{\ell}\mE\left[\|\nab \tilde{p}(t_{n+1})\|^2_{L^2}\right] &= \frac{Ch^2}{k}\sum_{n = 0}^{\ell}\mE\left[\|\nab(P(t_{n+1}) - P(t_n))\|^2_{L^2}\right]\leq Ch^2.
		\end{align*}
		
		Next, we estimate the second term $Z_2$. Using Lemma \ref{lemma2.2} and Lemma \ref{stability_pdes}, and also Lemma \ref{stability_FEMs}, we have
		\begin{align*}
			Z_2 &= \sum_{n = 0}^{\ell}\mE\biggl[C\int_{t_{n}}^{t_{n+1}} \left[\|\nab(\vu(t_{n+1}) - \vu(s))\|^2_{\vL^2} + \|\nab(\theta(t_{n+1}) - \theta(s))\|^2_{L^2}\right]\, ds  \\\nonumber
			&\qquad + \int_{t_{n}}^{t_{n+1}}\|\theta(t_{n+1}) - \theta(s)\|^2_{L^2}\, ds \\\nonumber
			&\qquad+ \frac{512C_e^4k}{\nu^2}\|\vu(t_{n+1}) - \vu(t_{n})\|^4_{\vH^1}\|\cQ_h\ve^{n+1}\|^2_{\vL^2} \\\nonumber
			&\qquad+  C\int_{t_n}^{t_{n+1}}\|\vu(t_{n}) - \vu(s)\|^2_{\vL^4}\|\vu(t_{n+1})\|^2_{\vL^4}\, ds \\\nonumber
			&\qquad+ C\int_{t_n}^{t_{n+1}}\|\vu(t_{n+1}) - \vu(s)\|^2_{\vL^4}\|\vu(s)\|^2_{\vL^4}\, ds
			\\\nonumber
			&\qquad+\frac{2C^2_ek}{\nu}\|\nab(\vu(t_{n+1})-\vu(t_n))\|^2_{\vL^2}\left[\|\cQ_h\ve^{n+1}\|^2_{\vL^2} + \|\cQ_h\ve^{n}\|^2_{\vL^2}\right] 
			\\\nonumber
			&\qquad+ C\int_{t_n}^{t_{n+1}}\|\vu(t_{n+1}) - \vu(s)\|^2_{\vH^1}\|\theta(t_{n+1})\|^2_{H^1}\, ds \\\nonumber
			&\qquad+ C\int_{t_n}^{t_{n+1}}\|\theta(t_{n+1}) - \theta(s)\|^2_{H^1}\|\vu(s)\|^2_{\vH^1}\, ds	
			\\\nonumber
			&\qquad+\frac{2C^2_ek}{\nu}\|\nab(\theta(t_{n+1})- \theta(t_{n}))\|^2_{L^2}\left[\|B_he_{\theta}^{n+1}\|^2_{L^2} + \|\cQ_h\ve^{n}\|^2_{\vL^2}\right]  \\\nonumber
			&\qquad+ \frac{512C_e^4k}{\nu^2}\|\theta(t_{n+1}) - \theta(t_{n})\|^4_{H^1}\|B_he_{\theta}^{n+1}\|^2_{L^2}\biggr]\\\nonumber
			&\leq Ck^{2\alpha} + C\left(\mE\left[\max_{1\leq n \leq M}\|\cQ_h\ve^{n}\|^4_{\vL^2}\|\vu(t_n)\|^4_{\vH^1}\right]\right)^{\frac12}\left(\mE\left[\left(k\sum_{n = 0}^{\ell} \|\vu(t_{n+1}) - \vu(t_n)\|^2_{\vH^1} \right)^{2}\right]\right)^{\frac12}\\\nonumber
			&\qquad+  C\left(\mE\left[\sup_{s\in [0,T]}\|\vu(s)\|^4_{\vL^4}\right]\right)^{\frac12}\left(\mE\left[\left(\sum_{n = 0}^{\ell} \int_{t_{n}}^{t_{n+1}}\|\vu(t_{n}) - \vu(s)\|^2_{\vL^4}\, ds\right)^{2}\right]\right)^{\frac12}\\\nonumber
			&\qquad+  C\left(\mE\left[\sup_{s\in [0,T]}\|\theta(s)\|^4_{\vH^1}\right]\right)^{\frac12}\left(\mE\left[\left(\sum_{n = 0}^{\ell} \int_{t_{n}}^{t_{n+1}}\|\vu(t_{n}) - \vu(s)\|^2_{\vH^1}\, ds\right)^{2}\right]\right)^{\frac12}\\\nonumber
			&\qquad+  C\left(\mE\left[\sup_{s\in [0,T]}\|\vu(s)\|^4_{\vH^1}\right]\right)^{\frac12}\left(\mE\left[\left(\sum_{n = 0}^{\ell} \int_{t_{n}}^{t_{n+1}}\|\theta(t_{n+1}) - \theta(s)\|^2_{H^1}\, ds\right)^{2}\right]\right)^{\frac12}\\\nonumber
			&\qquad+  C\left(\mE\left[\max_{1\leq n \leq M}\left(\|B_h e_{\theta}^n\|^4_{L^2} + \|\cQ_h\ve^n\|^4_{\vL^2}\right)\right]\right)^{\frac12}\left(\mE\left[\left(k\sum_{n = 0}^{\ell} \|\theta(t_{n+1}) - \theta(t_{n})\|^2_{H^1}\right)^{2}\right]\right)^{\frac12}\\\nonumber
			&\qquad+  C\left(\mE\left[\max_{1\leq n \leq M}\|B_h e_{\theta}^n\|^4_{L^2}\|\theta(t_n)\|^4_{H^1}\right]\right)^{\frac12}\left(\mE\left[\left(k\sum_{n = 0}^{\ell} \|\theta(t_{n+1}) - \theta(t_{n})\|^2_{H^1}\right)^{2}\right]\right)^{\frac12}\\\nonumber
			&\leq Ck^{2\alpha}.
		\end{align*}
		
		Finally, we estimate $Z_3$. Using the It\^o isometry, the martingale property of the It\^o integrals, Lemma \ref{lemma2.2}, the assumptions (B3), and (B4), and Lemma \ref{stability_pdes} we obtain
		
		\begin{align*}
			Z_3 &= \sum_{n = 0}^{\ell}\mE\biggl[2\left\|\int_{t_n}^{t_{n+1}}\left(\vG_1(\vu(s)) - \vG_1(\vu(t_n))\right)\, dW_1(s)\right\|^2_{\vL^2}\\\nonumber
			&\qquad+2\left\|\int_{t_n}^{t_{n+}} \left(G_2(\theta(s)) - G_2(\theta(t_n))\right)\, dW_2(s)\right\|^2_{L^2} 
			\\\nonumber
			&\qquad+  \mathbf{1}_{\Omega_{\rho, n}}\left(\int_{t_n}^{t_{n+1}}\left(\vG_1(\vu(s)) - \vG_1(\vu(t_n))\right)\, dW_1(s),  \cQ_h\ve^n\right)\\\nonumber
			&\qquad+  \mathbf{1}_{\Omega_{\rho, n}}\left(\int_{t_n}^{t_{n+}} \left(G_2(\theta(s)) - G_2(\theta(t_n))\right)\, dW_2(s),  B_h e^{n}_{\theta}\right)  \\\nonumber
			&\qquad+  \mathbf{1}_{\Omega_{\rho, n}}\left(\left(\vG_1(\vu(t_n))-\vG_1(\vu_h^n)\right)\Delta W_{1,n+1},\cQ_h\ve^{n}\right)  \\\nonumber
			&\qquad+  \mathbf{1}_{\Omega_{\rho, n}}\left(\left(G_2(\theta(t_n))-G_2(\theta_h^n)\right)\Delta W_{2,n+1}, B_h e_{\theta}^{n}\right) \\\nonumber
			&\qquad+ Ch^4 \left[\|\vA \vu(t_n)\|^2_{\vL^2}|\Delta W_{1,n+1}|^2 + \|\theta(t_n)\|^2_{H^2}|\Delta W_{2,n+1}|^2\right] \biggr]\\\nonumber
			&= \sum_{n = 0}^{\ell}\mE\biggl[2\int_{t_n}^{t_{n+1}}\left\|\vG_1(\vu(s)) - \vG_1(\vu(t_n))\right\|^2_{\vL^2}\, ds\\\nonumber
			&\qquad+2\int_{t_n}^{t_{n+}} \left\|G_2(\theta(s)) - G_2(\theta(t_n))\right\|^2_{L^2}\, ds\biggr] 
			\\\nonumber
			&\qquad+  0 + 0 + 0 + 0 + Ch^4\\\nonumber
			&\leq Ck^{2\alpha} + Ch^4.
		\end{align*}
		
		Now, we substitute all of the estimate from $Z_1$, $Z_2$ and $Z_3$ into the right-hand side of \eqref{eq_3.23} and then use the discrete Gronwall inequality to obtain
		\begin{align*}
			&	\frac12\mE\left[\mathbf{1}_{\Ome,\ell+1}\|\cQ_h\ve^{\ell+1}\|^2_{\vL^2} +  	\frac12\mathbf{1}_{\Ome,\ell+1}\|B_he_{\theta}^{\ell+1}\|^2_{L^2}\right] \\\nonumber
			&\qquad+ \frac{\nu k}{16}\sum_{n =0}^{\ell}\mE\left[\mathbf{1}_{\Omega_{\rho, n}}\|\nab\ve^{n+1}\|^2_{\vL^2} + \frac{\mu k}{2}\mathbf{1}_{\Omega_{\rho, n}}\|\nab e_{\theta}^{n+1}\|^2_{L^2}\right] \\\nonumber
			&\leq  Ck^{2\alpha} + Ch^2 + \sum_{n = 0}^{\ell}k \left(\frac32 + \tilde{C}\rho + 4C_{\vG_1}+ 4 C_{G_2}\right)\mE\left[\mathbf{1}_{\Omega_{\rho, n}}\left(\|\cQ_h\ve^{n+1}\|^2_{\vL^2} + \|\cQ_h\ve^{n}\|^2_{\vL^2} \right.\right.\\\nonumber
			&\qquad\qquad\qquad\left.\left.+ \|B_he_{\theta}^{n+1}\|^2_{L^2} + \|B_he_{\theta}^{n}\|^2_{L^2}\right)\right]\\\nonumber
			&\leq C\left(k^{2\alpha} + h^2 \right)e^{T\tilde{C}\rho}\\\nonumber
			&= C\ln(1/k)\left(k^{2\alpha} + h^{2}\right),
		\end{align*}
		where $\rho=  \frac{\ln(\ln(k^{-1}))}{\tilde{C}T}$  is used to obtain the last equality above.
		
	\end{proof}

	Finally, we state and prove the following error estimate for the pressure approximation:
	\begin{theorem}
		Under the same conditions of Theorem \ref{Main_Theorem_error}. There holds
		\begin{align}
			\left(\max_{1\leq n \leq M}\mE\left[\mathbf{1}_{\Omega_{\rho, n}}\left\|P(t_n) - k\sum_{\ell=1}^{n}p_h^{\ell}\right\|^2_{L^2}\right]\right)^{\frac12}  &\leq C\sqrt{\ln(1/k)}\left(k^{\alpha} + h\right),
		\end{align}
		where $C= C(\vu^0,\theta^0,T,\beta_1)$ is a postive constant.
	\end{theorem}
	\begin{proof}
		The proof is based on the inf-sup condition \eqref{inf-sup_discrete} and Theorem \ref{Main_Theorem_error}. First, denote $e^n_p: = P_h \tilde{p}(t_n) - p_h^n$. 
		Adding and subtracting the term $P_h \tilde{p}(t_n)$ in the error equations \eqref{eq_error_u} and then summing up the result from $n=0$ to $\ell$ for $1\leq \ell <M$, we obtain
		\begin{align*}
			\left(k\sum_{n=0}^{\ell} e^{n+1}_p, \div \pphi_h\right) &= \left(\ve^{\ell+1} - \ve^0, \pphi_h\right) + \nu \left(k\sum_{n = 0}^{\ell}\nab \ve^{n+1},\nab\pphi_h\right) \\\nonumber
			&\qquad- k\sum_{n = 0}^{\ell}\left(\tilde{p}(t_{n+1}) - P_h\tilde{p}(t_{n+1}), \div \pphi_h\right) \\\nonumber
			&\qquad-\sum_{n = 0}^{\ell} \biggl\{\int_{t_n}^{t_{n+1}} \nu\left(\nab(\vu(t_{n+1}) - \vu(s)),\nab \pphi_h\right)\, ds \\\nonumber
			&\qquad + k \left(e_{\theta}^n{\bf e}_2,\pphi_h\right) + \int_{t_{n}}^{t_{n+1}}\left((\theta(s) - \theta(t_n)){\bf e}_2,\pphi_h\right)\, ds \\\nonumber
			&\qquad +\int_{t_n}^{t_{n+1}}\left[{b}(\vu(t_{n}),\vu(t_{n+1}),\pphi_h) - {b}(\vu(s), \vu(s),\pphi_h) \right] \, ds
			\\\nonumber
			&\qquad -k\left[{b}(\vu(t_{n}),\vu(t_{n+1}),\pphi_h) - {b}(\vu_h^{n}, \vu_h^{n+1},\pphi_h)  \right]\biggr\}\\\nonumber
			&\qquad -\sum_{n = 0}^{\ell}\biggl\{\left(\int_{t_n}^{t_{n+1}}\left(\vG_1(\vu(s)) - \vG_1(\vu(t_n))\right)\, dW_1(s), \pphi_h\right)\\\nonumber
			&\qquad+\left(\left(\vG_1(\vu(t_n))-\vG_1(\vu_h^n)\right)\Delta W_{1,n+1},\pphi_h\right)\biggr\}.
		\end{align*}
		
		Next, using the inf-sup condition \eqref{inf-sup_discrete}, and Cauchy-Schwarz's inequality we obtain
		\begin{align*}
			\beta_1 \left\|k\sum_{n=0}^{\ell} e^{n+1}_p\right\|_{L^2} &\leq \sup_{\pphi_h \in \mH_h, \pphi_h \neq 0}\frac{\left(k\sum_{n=0}^{\ell} e^{n+1}_p, \div \pphi_h\right) }{\|\pphi_h\|_{\mH_h}}\\\nonumber
			&\leq C \|\ve^{\ell+1} - \ve^0\|_{\vL^2} + \nu \left\|k\sum_{n = 0}^{\ell} \nab \ve^{n+1}\right\|_{\vL^2} \\\nonumber
			&\qquad+ Ck\sum_{n = 0}^{\ell} \|\tilde{p}(t_{n+1}) - P_h\tilde{p}(t_{n+1})\|_{L^2}\\\nonumber
			&\qquad+ C\sum_{n =0}^{\ell}\frac{1}{\|\pphi_h\|_{\mH_h}}\biggl\{\int_{t_n}^{t_{n+1}} \nu\left(\nab(\vu(t_{n+1}) - \vu(s)),\nab \pphi_h\right)\, ds \\\nonumber
			&\qquad + k \left(e_{\theta}^n{\bf e}_2,\pphi_h\right) + \int_{t_{n}}^{t_{n+1}}\left((\theta(s) - \theta(t_n)){\bf e}_2,\pphi_h\right)\, ds \\\nonumber
			&\qquad +\int_{t_n}^{t_{n+1}}\left[{b}(\vu(t_{n}),\vu(t_{n+1}),\pphi_h) - {b}(\vu(s), \vu(s),\pphi_h) \right] \, ds
			\\\nonumber
			&\qquad -k\left[{b}(\vu(t_{n}),\vu(t_{n+1}),\pphi_h) - {b}(\vu_h^{n}, \vu_h^{n+1},\pphi_h)  \right]\biggr\}\\\nonumber
			&\qquad + C\sum_{n = 0}^{\ell}\frac{1}{\|\pphi_h\|_{\mH_h}}\biggl\{\left(\int_{t_n}^{t_{n+1}}\left(\vG_1(\vu(s)) - \vG_1(\vu(t_n))\right)\, dW_1(s), \pphi_h\right)\\\nonumber
			&\qquad+\left(\left(\vG_1(\vu(t_n))-\vG_1(\vu_h^n)\right)\Delta W_{1,n+1},\pphi_h\right)\biggr\}\\\nonumber
			&:=I + II + III + IV+ V.
		\end{align*}
		
		Using the Cauchy-Schwarz inequality, Lemma \ref{lemma2.2} and Theorem \ref{Main_Theorem_error}, we easily obtain the following estimates
		\begin{align}\label{eq_3.22}
			\left\{\mE\left[\mathbf{1}_{\Omega_{\rho,\ell}}(I^2 + II^2 + IV^2 + V^2)\right]\right\}^{1/2}\leq C\sqrt{\ln(1/k)}\left(k^{\alpha} + h\right).
		\end{align}
		
		Next, using \eqref{P_property_a} and Lemma \ref{lemma2.2} (d), we have
		\begin{align}\label{eq_3.23}
			\left\{\mE[III^2]\right\}^{1/2} &\leq C\left\{k\sum_{n = 0}^{\ell} \mE\left[\|\tilde{p}(t_{n+1}) - P_h\tilde{p}(t_{n+1})\|^2_{L^2}\right]\right\}^{1/2} \\\nonumber
			&\leq Ch \left\{k\sum_{n = 0}^{\ell} \mE[\|\nab \tilde{p}(t_{n+1})\|^2_{L^2}]\right\}^{1/2} \\\nonumber
			&= Ch \left\{k^{-1}\sum_{n = 0}^{\ell} \mE[\|\nab(P(t_{n+1}) - P(t_n))\|^2_{L^2}]\right\}^{1/2} \leq Ch.
		\end{align}
		
		Combining both \eqref{eq_3.22} and \eqref{eq_3.23}, we have
		\begin{align}\label{eq_3.24}
			\left(\mE\left[\mathbf{1}_{\Omega_{\rho,\ell}}\left\|k\sum_{\ell=1}^{\ell} e_p^{n}\right\|^2_{L^2}\right]\right)^{1/2} \leq C\sqrt{\ln(1/k)}\left(k^{\alpha} + h\right).
		\end{align}
		
		Finally, the proof is completed by using the triangle inequality with the decomposition
		\begin{align*}
			P(t_{\ell}) - k\sum_{n = 1}^{\ell} p_{h}^n= k\sum_{n = 1}^{\ell}\left(\tilde{p}(t_n) - p_h^n\right) = k\sum_{n = 1}^{\ell}\left(\tilde{p}(t_n) - P_h\tilde{p}(t_n)\right) + k\sum_{n = 1}^{\ell} e^n_p,
		\end{align*}
		and \eqref{eq_3.24}, \eqref{P_property}, and \eqref{P_property_a}.
	\end{proof}	
	
	\section{Numerical experiments}	\label{section 4}
	
	In this section, we present numerical experiments to validate our theoretical results. Throughout all tests, we set $D = (0,1)^2 \subset \mathbb{R}^2$, $T = 1$, $\nu = 1$, and $\mu = 1$. 
	The Wiener processes $W_1(t)$ and $W_2(t)$ in \eqref{eq1.1} are taken to be $\mathbb{R}$-valued and are simulated using a minimal time-step size $k_0 = 2^{-11}$. We also choose $\vG_1(\vu) = \vu$ and $G_2(\theta) = \theta$ in all the numerical tests.
	For all computations, we approximate expectations using the standard Monte Carlo method with $J = 400$ independent samples. 
	Spatial discretization is performed using the MINI element for the mixed element pair of $(\vu_h^n, p_h^n)$, while the linear polynomial is chosen for the element space of $\theta_h^n$.
	Finally, homogeneous Dirichlet boundary conditions are imposed on both $\vu$ and $\theta$ in all numerical tests.

	We implement the Main Algorithm and compute the errors of the velocity, pressure, and temperature approximations in the specified norms below. Since the exact solutions are unknown, the errors are computed between the computed solution $(\vu_h^n(\omega_j), p_h^n(\omega_j), \theta_h^n(\omega_j))$ and a reference solution $(\vu^n_{ref}(\omega_j), {{p}}^n_{ref}(\omega_j), {{\theta}}^n_{ref}(\omega_j))$ (specified later) at the $\omega_j$-th sample.
	
	Furthermore, to evaluate errors in strong norms, we use the following numerical 
	integration formulas: 
	\begin{align*}
		L^2_{\omega}L^{\infty}_tL^2_x(\vu)&:=	\bigl(\mE\bigl[\max_{1\leq n \leq M}\|\vu(t_n) - \vu^n_h\|^2_{\vL^2}\bigr]\bigr)^{1/2} \\
		&\approx \Bigl(\frac{1}{J}\sum_{j= 1}^J\bigl(\max_{1\leq n \leq M} \|\vu^n_{ref}(\omega_j) - \vu^n_h(\omega_j)\|^2_{\vL^2}\bigr)\Bigr)^{1/2},\\
		L^q_{\omega}L^{\infty}_tL^2_x(\theta)&:=	\bigl(\mE\bigl[\max_{1\leq n \leq M}\|\theta(t_n) - \theta^n_h\|^2_{L^2}\bigr]\bigr)^{1/2} \\
		&\approx \Bigl(\frac{1}{J}\sum_{j= 1}^J\bigl(\max_{1\leq n \leq M} \|\theta^n_{ref}(\omega_j) - \theta^n_h(\omega_j)\|^2_{L^2}\bigr)\Bigr)^{1/2},\\
		L^2_{\omega}L^2_tH^1_x(\vu)&:= \Bigl(\mE\Bigl[k\sum_{n = 1}^M \|\vu(t_n) - \vu^n_h\|^2_{\vH^1}\Bigr]\Bigr)^{1/2} \\
		&\approx \Bigl(\frac{1}{J}\sum_{j= 1}^J\Bigl(k\sum_{n = 1}^M\|\vu^n_{ref}(\omega_j) - \vu_h^n(\omega_j)\|^2_{\vH^1}\Bigr)\Bigr)^{1/2},\\
		L^2_{\omega}L^2_tH^1_x(\theta)&:= \Bigl(\mE\Bigl[k\sum_{n = 1}^M \|\theta(t_n) - \theta^n_h\|^2_{H^1}\Bigr]\Bigr)^{1/2} \\
		&\approx \Bigl(\frac{1}{J}\sum_{j= 1}^J\Bigl(k\sum_{n = 1}^M\|\theta^n_{ref}(\omega_j) - \theta_h^n(\omega_j)\|^2_{H^1}\Bigr)\Bigr)^{1/2},\\
		L^2_{\omega}L^2_xL^1_t(p) &: = \Bigl(\mE\Bigl[\Bigl\|P(t_M) -k\sum_{n = 1}^Mp^n_h \Bigr\|^2_{L^2}\Bigr]\Bigr)^{1/2}\\
		&\approx \Bigl( \frac{1}{J}\sum_{j= 1}^J \Bigl( \Bigl\|k\sum_{n = 1}^M\bigl(p^n_{ref}(\omega_j) - p^n_h(\omega_j)\bigr)\Bigr\|^2_{L^2}\Bigr)\Bigr)^{1/2}.
	\end{align*}
	
	\textbf{Test 1.}
	In this experiment, we investigate the temporal convergence rate \(O(k^{\alpha})\) for \(\alpha \in (0,1/2)\) predicted by Theorem~\ref{Main_Theorem_error}. To this end, we apply the Main Algorithm to compute numerical approximations \(\{(u_h^n, p_h^n, \theta_h^n)\}\) using a fixed spatial mesh size \(h = 1.67\times 10^{-2}\), while systematically refining the time step according to \(k = 2^{\ell} k_0\) for \(\ell \in \mathbb{N}\).
	
	The reference solutions \(\{(\vu_{\mathrm{ref}}^n, p_{\mathrm{ref}}^n, \theta_{\mathrm{ref}}^n)\}\) are obtained using a finer time step \(k_{\mathrm{ref}} = k/2\). The temporal discretization errors are then approximated by comparing numerical solutions computed at two successive time resolutions. This approach is standard in numerical studies of stochastic Navier--Stokes equations and related systems; see, for example, \cite{vo2025high, FV2024, vo2023higher}.
	
\begin{figure}[hbt]
		\begin{center}
			\includegraphics[scale=0.11]{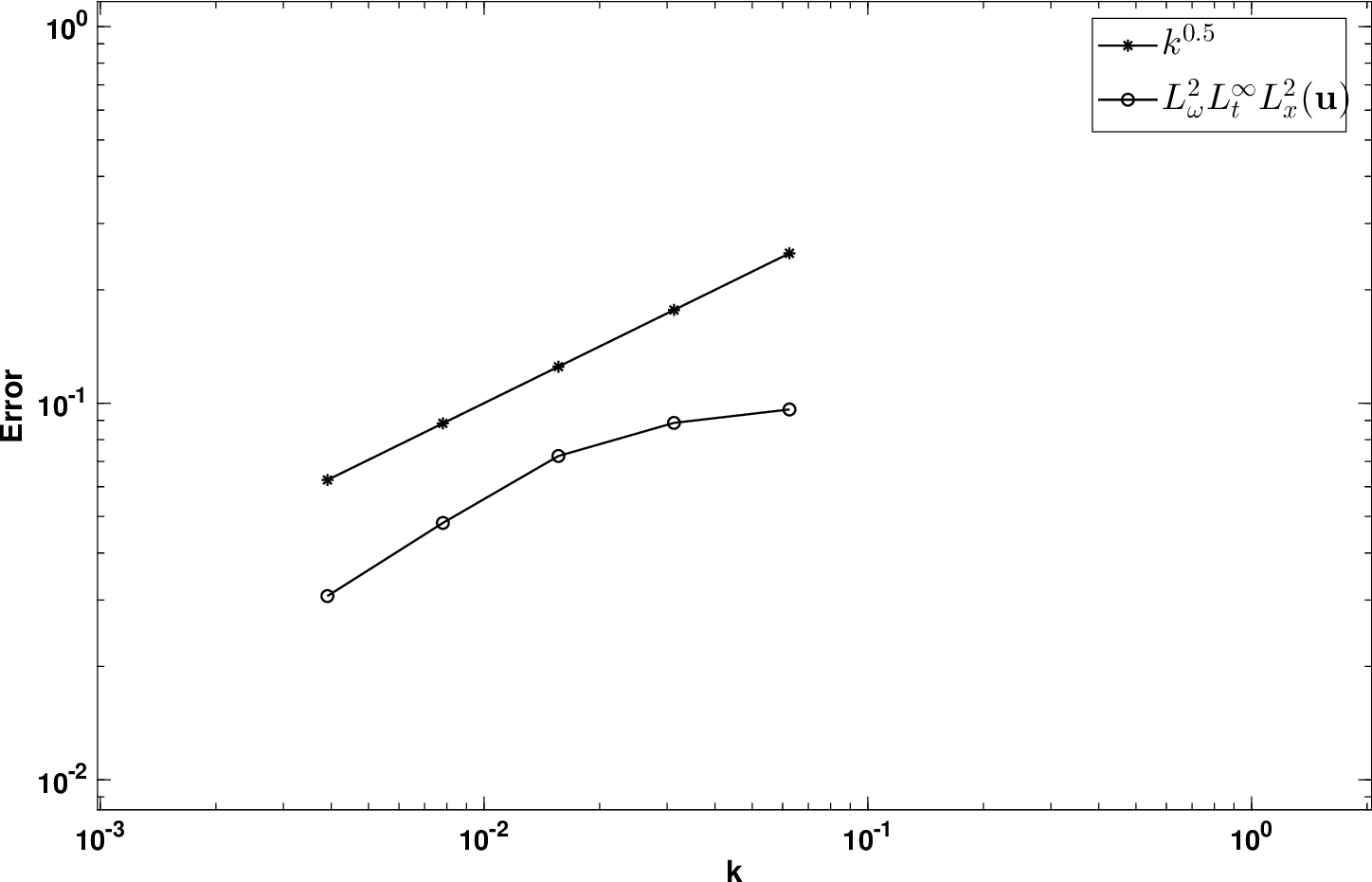}
			\includegraphics[scale=0.11]{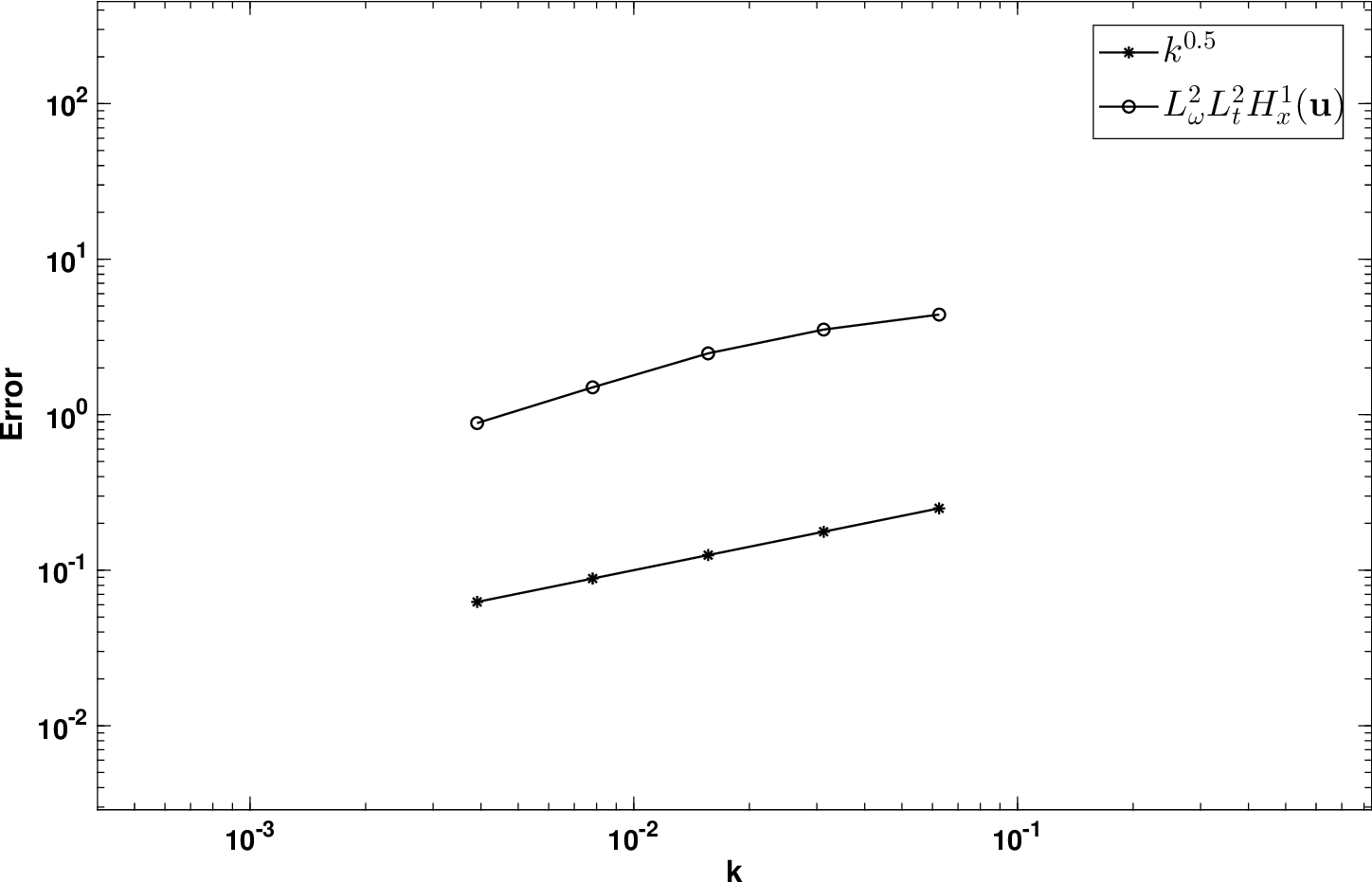}
			\caption{Plots of the time discretization errors and convergence order of the computed velocity $\{ {\bf u}^n_h\}$ in $L^2$-norm(left) and $H^1$-norm (right) with $k = 2^{-4}, 2^{-5}, 2^{-6}, 2^{-7}, 2^{-8} $.}\label{fig4.1}
		\end{center}
	\end{figure}

\begin{figure}[hbt]
		\begin{center}
			\includegraphics[scale=0.11]{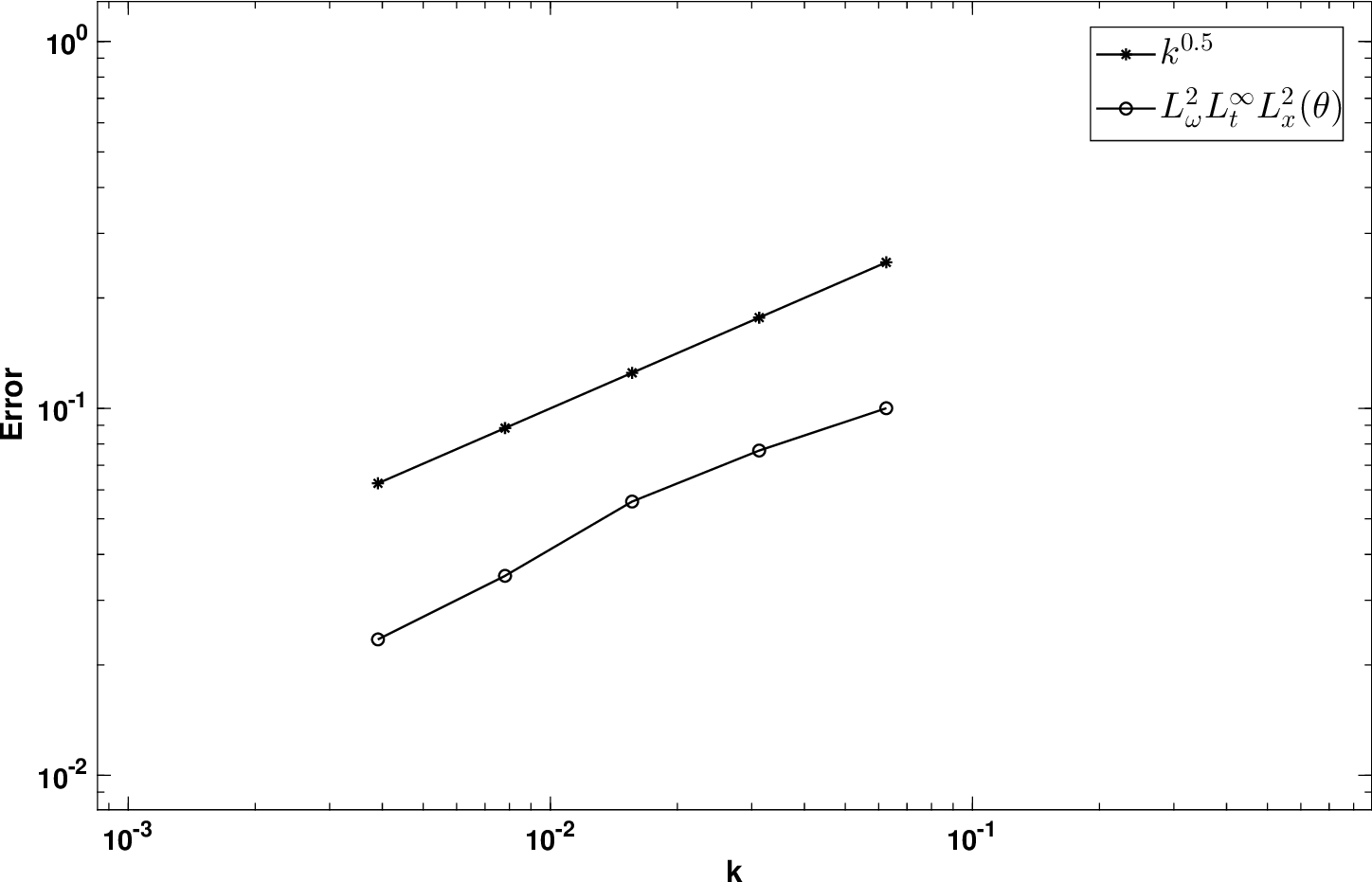}
			\includegraphics[scale=0.11]{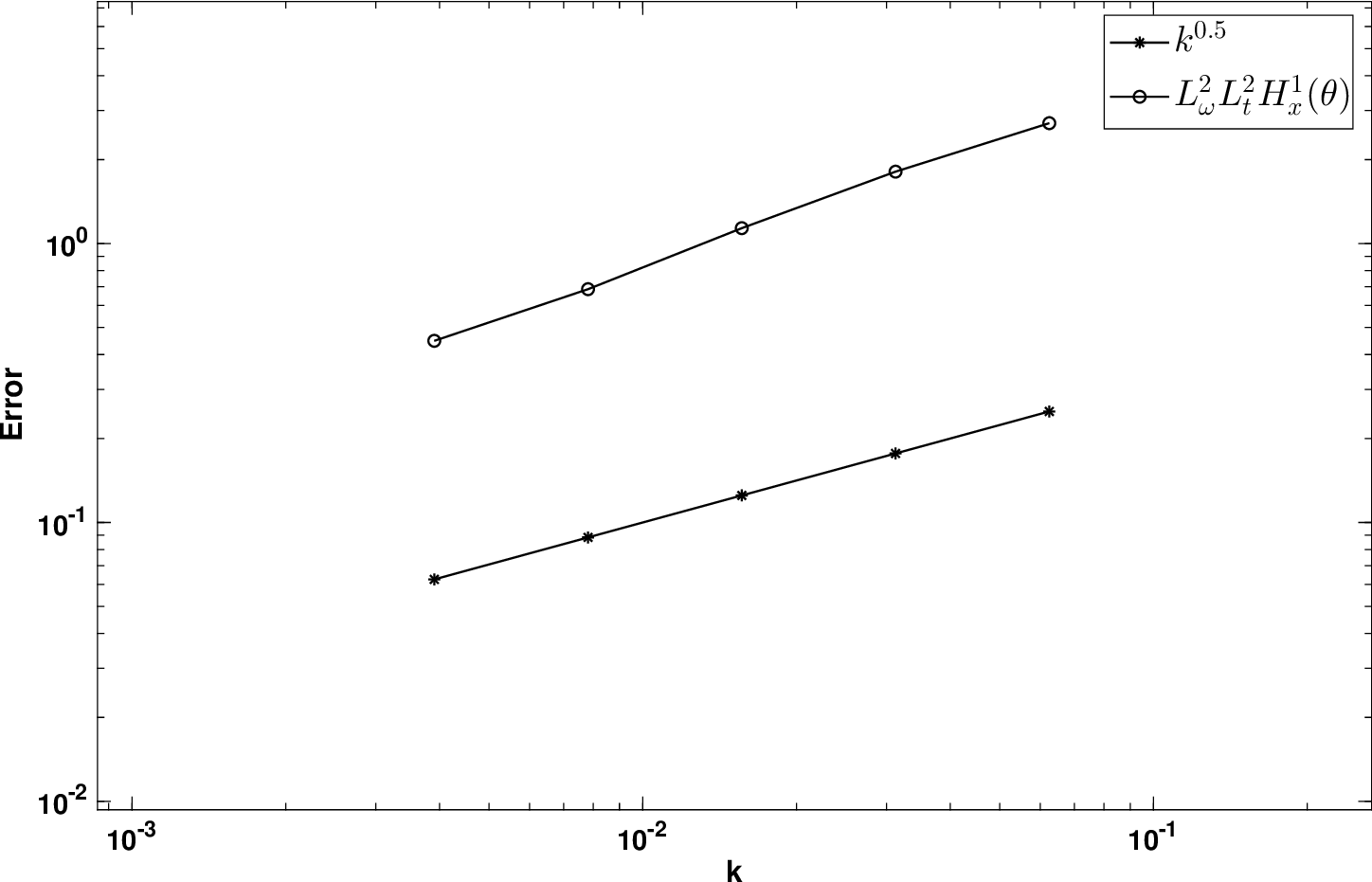}
			\caption{Plots of the time discretization errors and convergence order of the computed temperature $\{ {\theta}^n_h\}$ in $L^2$-norm(left) and $H^1$-norm (right) with $k = 2^{-4}, 2^{-5}, 2^{-6}, 2^{-7}, 2^{-8} $.}\label{fig4.2}
		\end{center}
	\end{figure}

\begin{figure}[hbt]
		\begin{center}
			\includegraphics[scale=0.11]{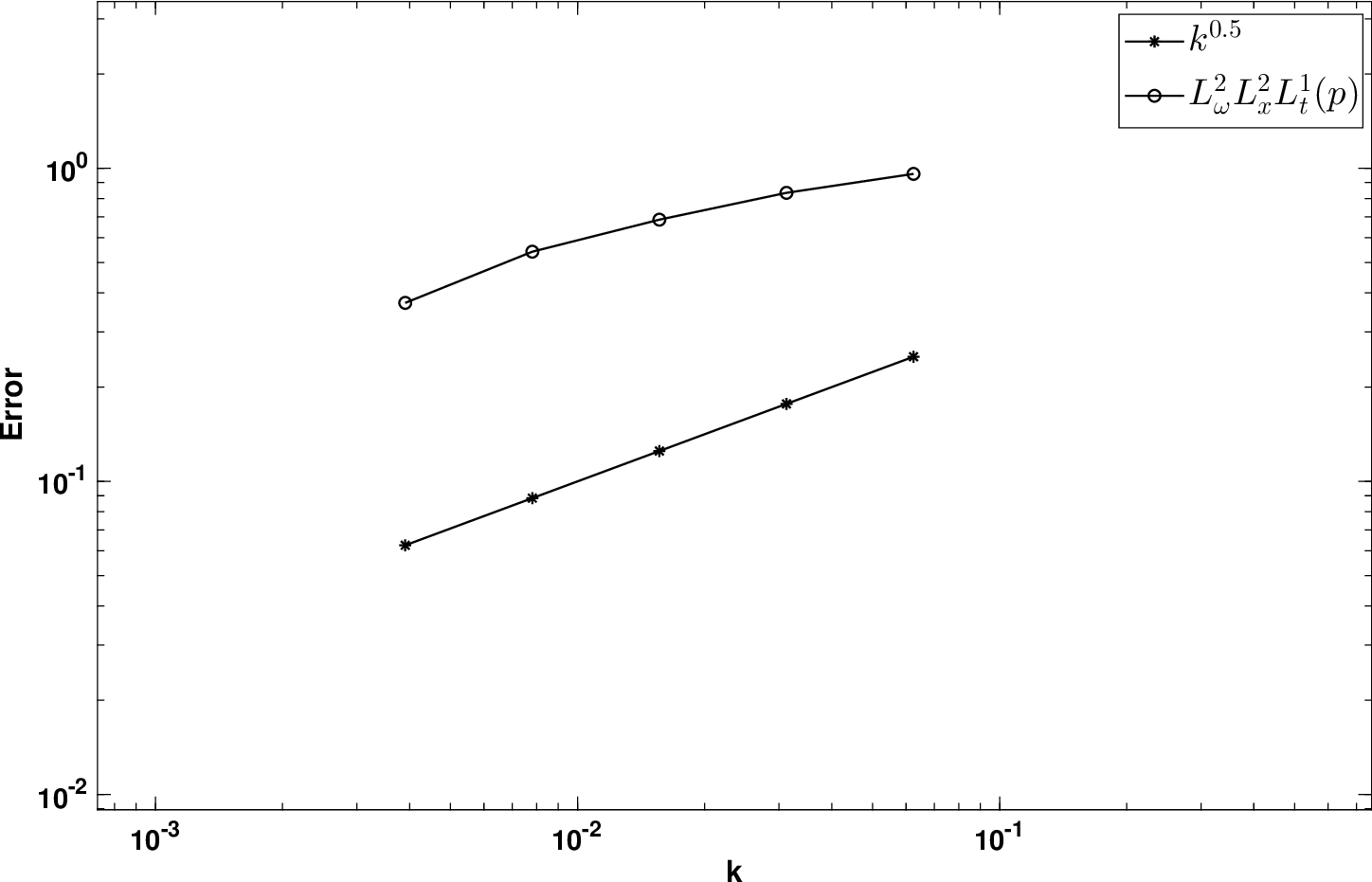}
			\caption{Plots of the time discretization errors and convergence order of the computed pressure $\{ {\bf u}^n_h\}$ in time-averaged norm with $k = 2^{-4}, 2^{-5}, 2^{-6}, 2^{-7}, 2^{-8} $.}\label{fig4.3}
		\end{center}
	\end{figure}

	The resulting errors are shown in Figures~\ref{fig4.1}--\ref{fig4.3}. The observed convergence rates for the velocity, pressure, and temperature are consistent with the theoretical prediction of an order close to \(1/2\) stated in Theorem~\ref{Main_Theorem_error}.

	
	\medskip

	\textbf{Test 2.} Finally, we present a numerical example to validate the spatial convergence rates of the proposed numerical scheme predicted by Theorem~\ref{Main_Theorem_error}. In this experiment, the time step is fixed at $k = 2^{-9}$, and the spatial mesh size $h$ is successively refined.
	
\begin{figure}[hbt]
		\begin{center}
			\includegraphics[scale=0.11]{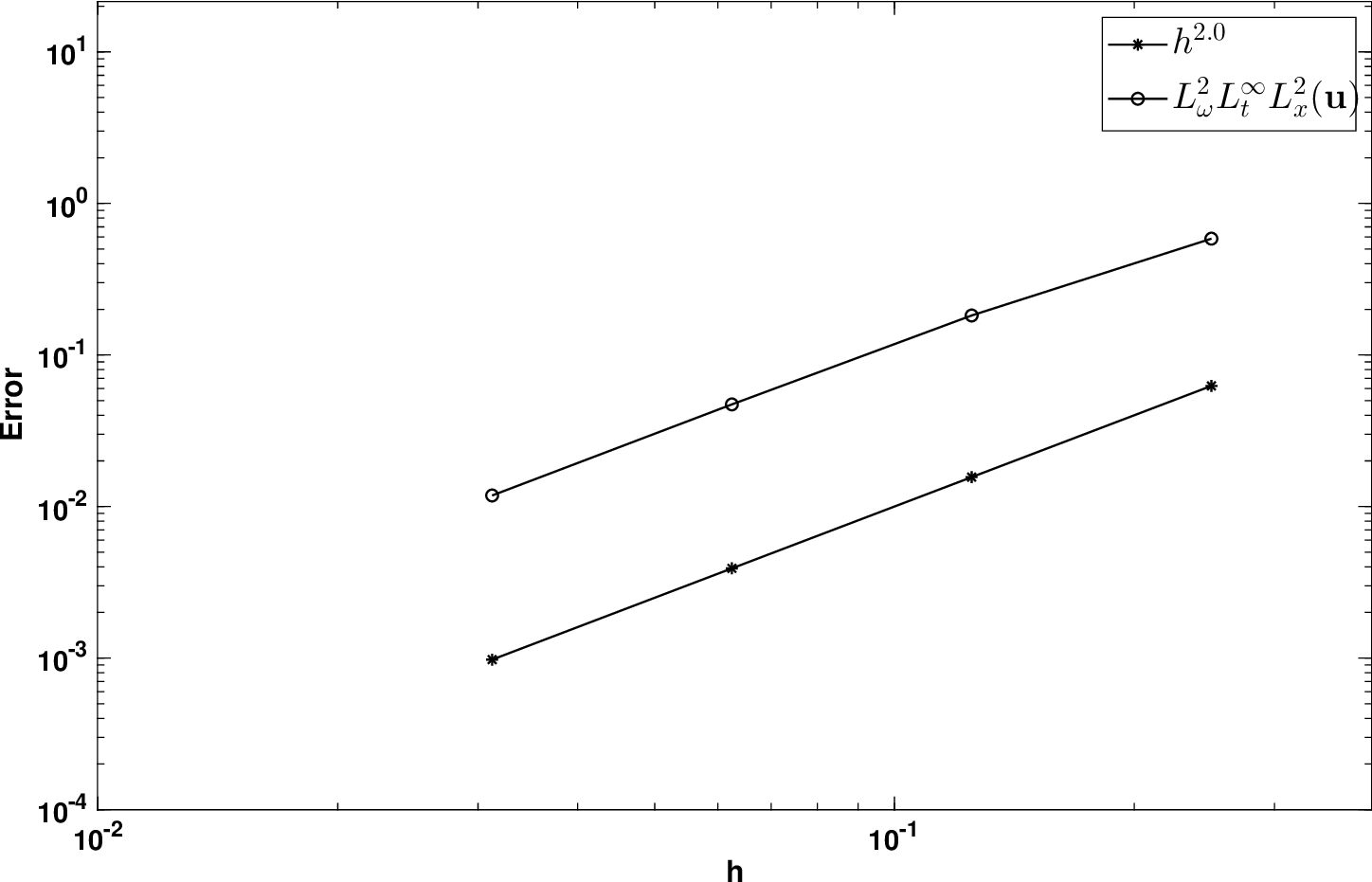}
			\includegraphics[scale=0.11]{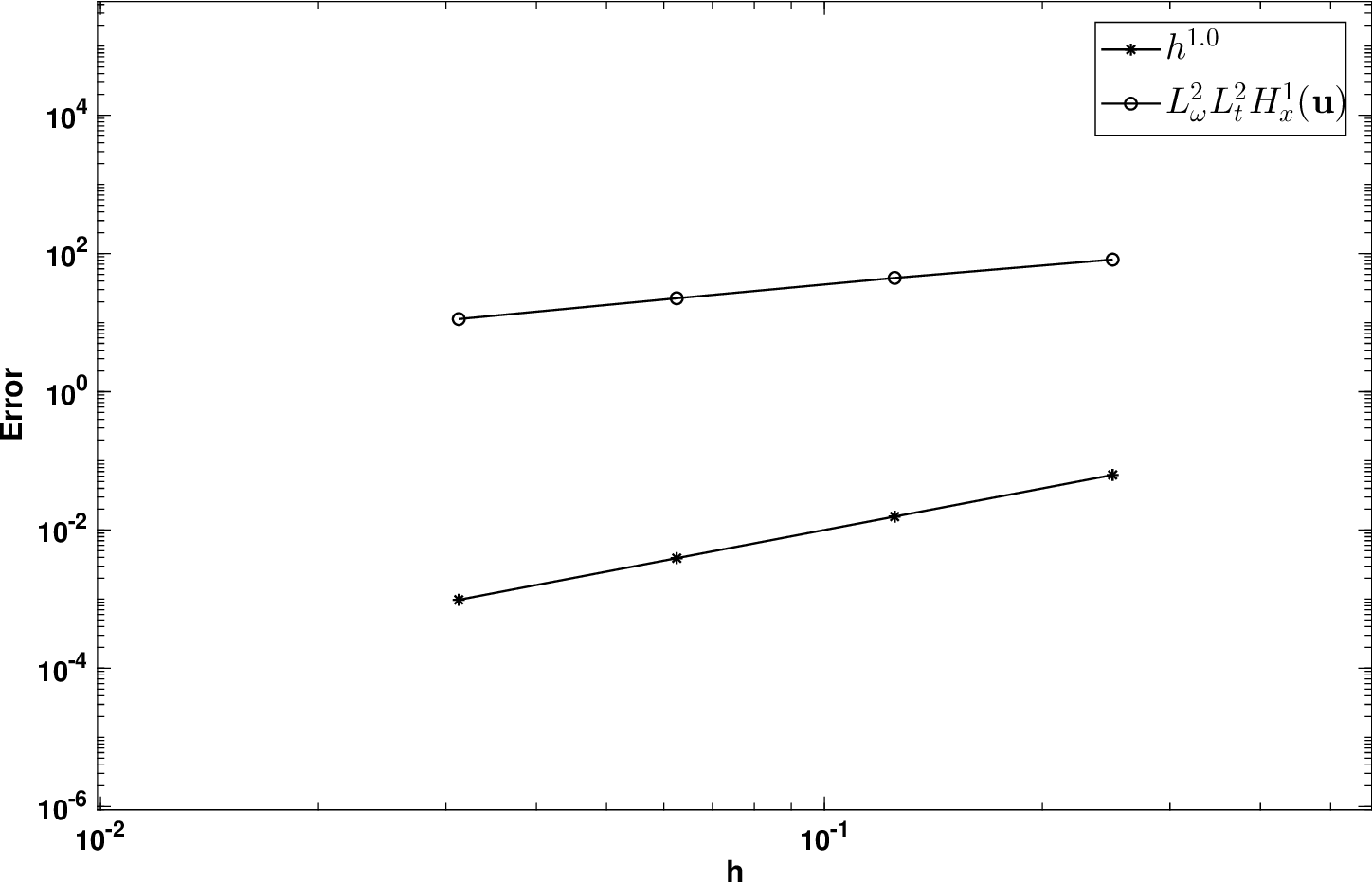}
			\caption{Plots of the spatial discretization errors and convergence order of the computed velocity $\{ {\bf u}^n_h\}$ in $L^2$-norm(left) and $H^1$-norm (right) with $h = 2^{-2}, 2^{-3}, 2^{-4}, 2^{-5}$.}\label{fig4.4}
		\end{center}
	\end{figure}

\begin{figure}[hbt]
		\begin{center}
			\includegraphics[scale=0.11]{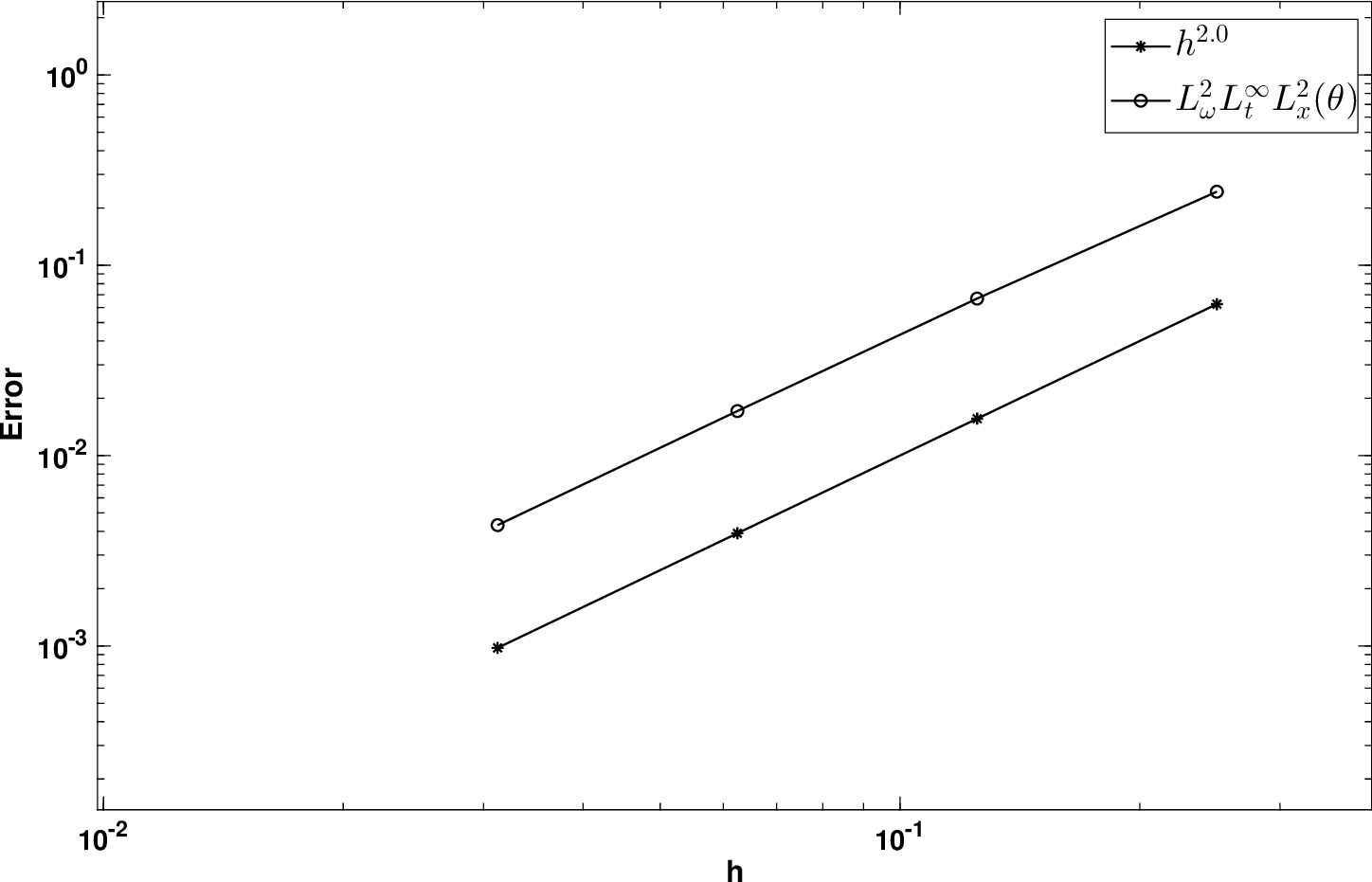}
			\includegraphics[scale=0.11]{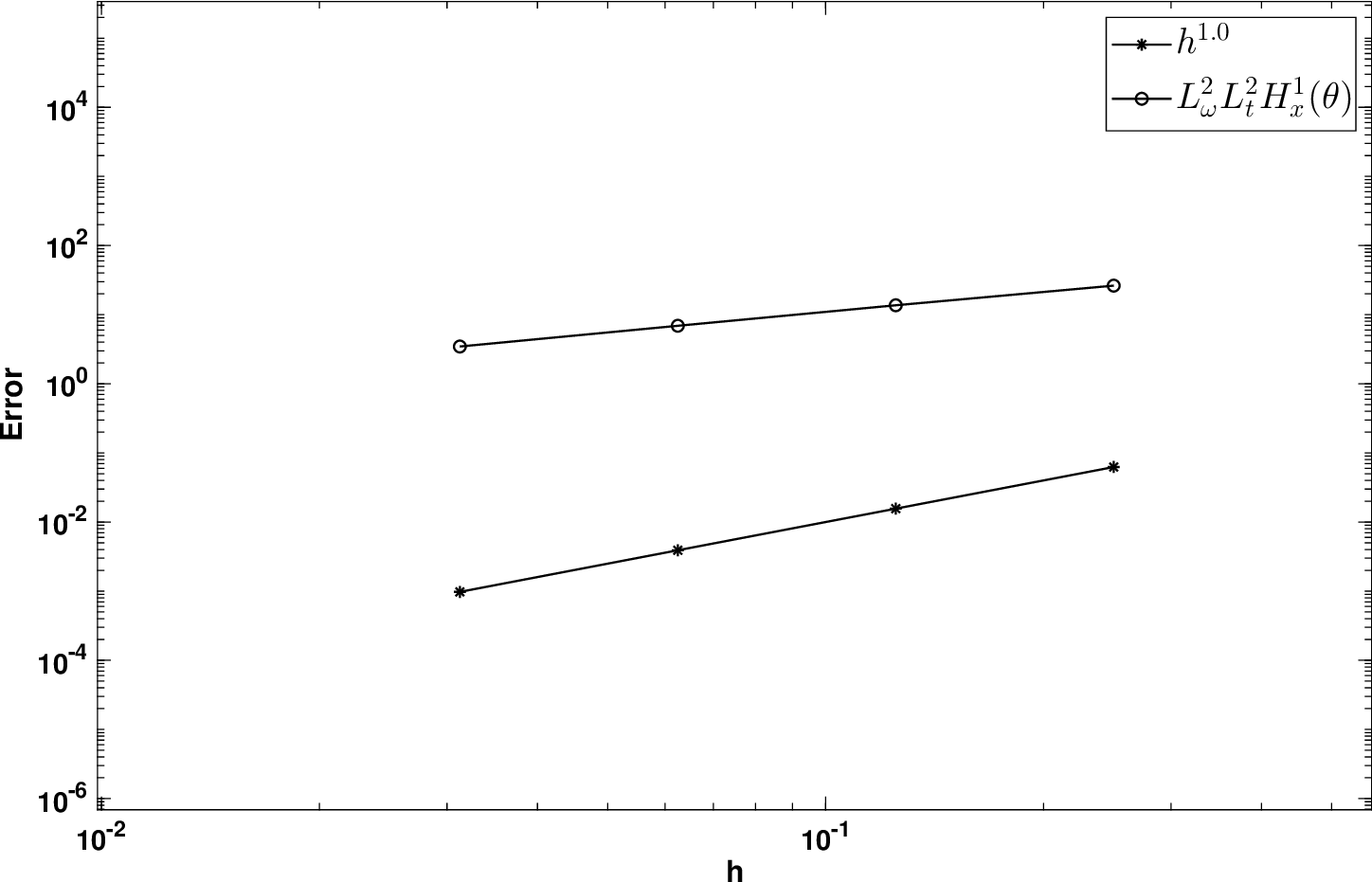}
			\caption{Plots of the spatial discretization errors and convergence order of the computed temperature $\{ {\theta}^n_h\}$ in $L^2$-norm(left) and $H^1$-norm (right) with $h = 2^{-2}, 2^{-3}, 2^{-4}, 2^{-5}$.}\label{fig4.5}
		\end{center}
	\end{figure}

\begin{figure}[hbt]
		\begin{center}
			\includegraphics[scale=0.11]{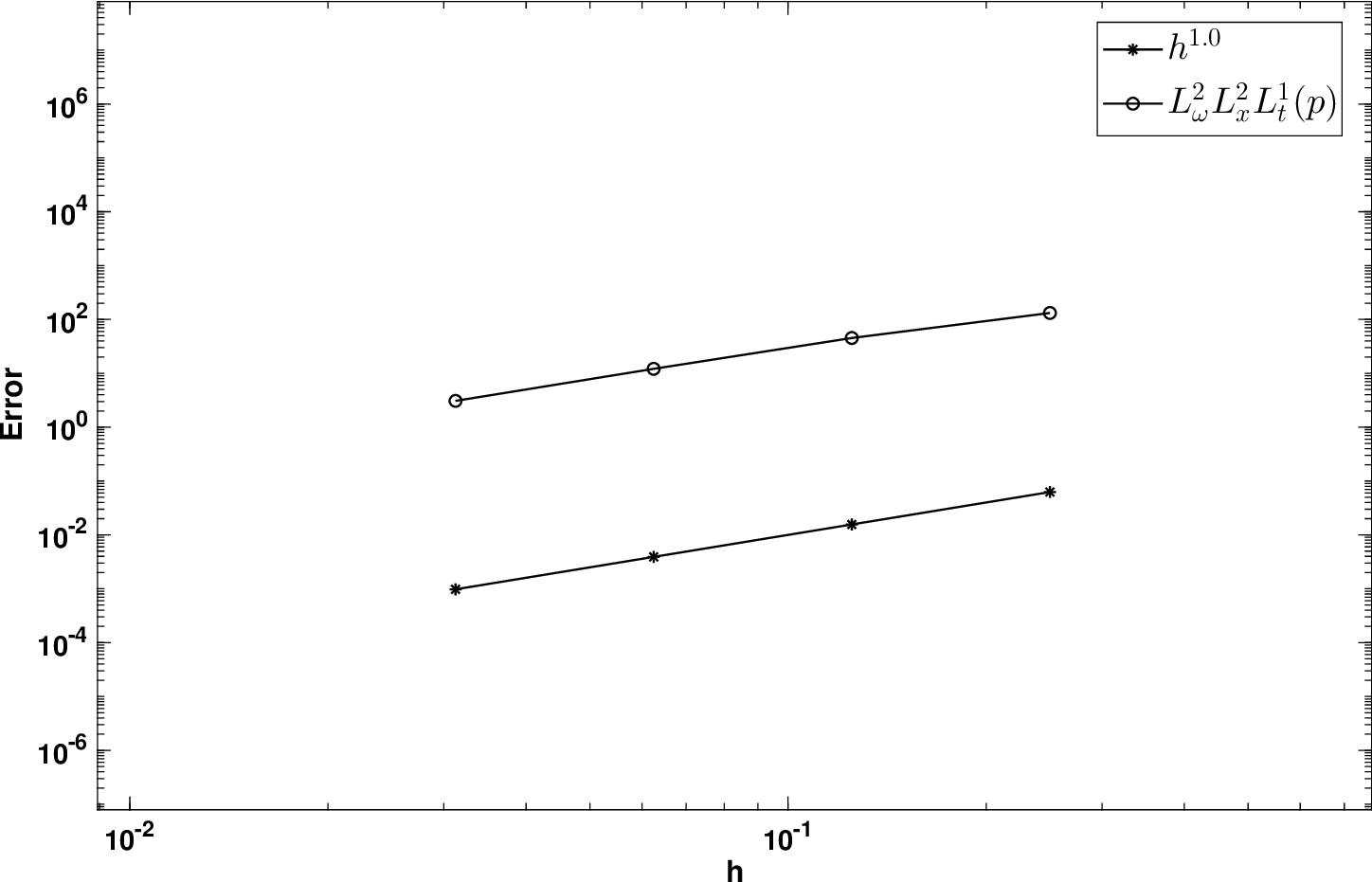}
			\caption{Plots of the spatial discretization errors and convergence order of the computed pressure $\{ {\bf u}^n_h\}$ in time-averaged norm with $h = 2^{-2}, 2^{-3}, 2^{-4}, 2^{-5}$.}\label{fig4.6}
		\end{center}
	\end{figure}

	The numerical results shown in Figures~\ref{fig4.4} (right) and \ref{fig4.5} (right) demonstrate a first-order convergence rate for the velocity approximation in the energy norm, in agreement with the theoretical result of Theorem~\ref{Main_Theorem_error}. Furthermore, Figures~\ref{fig4.4} (left) and \ref{fig4.5} (left) indicate a second-order convergence rate for the velocity approximation in the $L^2$ norm, which exceeds the theoretical estimate. Finally, Figure~\ref{fig4.6} confirms a first-order convergence rate for the pressure approximation.


	\section*{Declarations of Funding}\,
	
	The author Liet Vo was supported by the National Science Foundation (NSF) under Grant No. DMS-2530211.
	
	\section*{Data Availability} \,
	
Data sharing not applicable to this article as no datasets were generated or analysed during the current study.

\section*{Declaration of Conflict of Interest}\,
	
	The author has no conflict of interest.

	\bibliographystyle{abbrv}
	\bibliography{references_Boussinesq}

\end{document}